\documentclass{elsart}%
\usepackage{amsmath}
\usepackage{amsfonts}
\usepackage{amssymb}
\usepackage{graphicx}%
\newcommand{\citet}[1]{\cite{#1}}
\begin{document}
\begin{frontmatter}
\title{A wildland fire model with data assimilation}
\date{June 2006, revised January 2008.}
\author{Jan Mandel\corauthref{cor1}\thanksref{ucd}\thanksref{ncar}},
\author{Lynn S. Bennethum\thanksref{ucd}},
\author{Jonathan D. Beezley\thanksref{ucd}},
\author[ncar]{Janice L. Coen},
\author{Craig C.\ Douglas\thanksref{uky}\thanksref{yale}},
\author{Minjeong Kim\thanksref{ucd}},
and \author[av]{Anthony Vodacek}
\corauth[cor1]{Corresponding author. Campus Box 170, Denver CO 80217-3364, USA}
\address[ucd]{Center of Computational Mathematics and
Department of Mathematical Sciences, University of Colorado Denver, Denver, CO}
\address[ncar]{Mesoscale and Microscale Meteorology Division,\\ National Center for Atmospheric Research, Boulder, CO}
\address[uky]{Department of Computer Science, University of Kentucky,
Lexington, KY}
\address[yale]{Department of Computer Science, Yale University, New Haven, CT}
\address[av]{Center for Imaging Science, Rochester Institute
of Technology, Rochester, NY}
\begin{abstract}
A wildfire model is formulated based on balance equations for energy and fuel,
where the fuel loss due to combustion corresponds to the fuel reaction rate.
The resulting coupled partial differential equations have coefficients that
can be approximated from prior measurements of wildfires. An ensemble Kalman
filter technique with regularization is then used to assimilate temperatures measured at selected
points into running wildfire simulations. The assimilation technique is able
to modify the simulations to track the measurements correctly even if the
simulations were started with an erroneous ignition location that is quite far
away from the correct one.
\end{abstract}
\end{frontmatter}

%automatically added by CVS, do not edit!!

%\begin{verbatim}
%$Revision $Id: fire_simple.tex,v 1.235 2008/01/15 23:00:49 jmandel Exp $
%\end{verbatim}

\section{Introduction}

\label{sec:introduction}

Modeling forest fires is a multi-scale multi-physics problem. One can try to
account for the many physical processes involved at the appropriate scales,
but at the cost of speed. Simplifying appropriate physical processes is one
way to obtain a faster-running model. In this paper we also propose using a
data assimilation technique to incorporate data in real time. The purpose of
this paper is a demonstration-of-concept: we take a very simple model, develop
a data assimilation technique, and show how, even with this very simple model,
realistic results can be obtained even with significant errors in the initial
conditions (location of the fire). This is the first step of a longer-term
goal in which a more realistic model will be used.

Data assimilation is a technique used to incorporate data into a running model
using sequential statistical estimation. Data assimilation is made necessary
by the facts that no model is perfect, the available data is spread over time
and space, and it is burdened with errors. The model solution is just one
realization from a probability distribution. Data
assimilation methods have achieved good success in fields like oil and gas
pipeline models \cite{Emara-2002-NMS} and atmospheric, climate, and ocean
modeling \cite{Kalnay-AMD-2003}, and they are a part of virtually any
navigation system, from steering the Apollo moon spaceships in the 60's to GPS
and operating unmanned drones or rovers in hostile conditions like Afghanistan
or Mars today. Data assimilation can also dynamically steer the measurement
process, by suggesting locations where the collection of data would result in
the best reduction of uncertainty in the forecast~\cite{Kalnay-AMD-2003}.

A new paradigm in modeling beyond current techniques in data assimilation is
to use Dynamic Data-Driven Application System (DDDAS) techniques
\cite{Darema-2004-DDD}. Data assimilation is just one of the techniques from
the DDDAS toolbox, which entails the ability to dynamically incorporate
additional data into an executing application, and in reverse, the ability of
an application to dynamically steer the measurement process. Other DDDAS
techniques include deterministic methods such as time rollback, checkpointing,
data flow computations, and optimization. One aspect of DDDAS is using data
assimilation and measurement steering techniques from weather forecasting in
other fields. In a DDDAS, simulations and measurements become a symbiotic
feedback control system. Such capabilities promise more accurate analysis and
prediction, more precise controls, and more reliable outcomes.

Our ultimate objective is to build a real-time coupled atmospheric-wildland
fire modeling system based on DDDAS techniques that is steered dynamically by
data, where data includes atmospheric, fire, fuel, terrain, and other data
that influence the growth of fires
\cite{Douglas-2006-DVW,Mandel-2007-DDD,Mandel-2005-TDD,Mandel-2004-NDD}. This
work describes one stage of our investigation, that is, to develop and
validate techniques to ingest fire data that might originate from in situ and
remote sensors into a newly developed fire model. The purpose of this paper is
to combine a data assimilation method with a partial differential equation
(PDE) based model for wildland fire that achieves good physical behavior and
can run faster than real time. The model in this paper does not yet include
coupling with the atmosphere, though it is known that such coupling is
essential for the wildland fire behavior \cite{Clark-1996-CAF}. Coupling the
fire model with atmospheric dynamics as well as with data assimilation is
currently under development. Models using explicit, detailed combustion
physics are not feasible for prediction, since they require a large number of
chemical reactions and species and extremely high resolution (grid cells
$<<1m$) fluid dynamics \cite{Zhou-2001-ERM}. The actual interaction between
the atmosphere and the fire and vegetation is quite complicated,
involving turbulence in the vegetation layer and its consequences on heat
transfer and combustion \cite{Baines-1990-PMP}.
%the following remark copied verbatim from referee's report
The example provided in this paper is for $250\times250$ cells of $2m$ size
each. A model like FIRETEC \cite{Linn-2002-SWB,Linn-1996-FTM} could do the
same, but including the full interaction between fire, vegetation, and the
atmosphere, and it would come at a much higher computational cost. FIRETEC is an
example of a physically-based model that simplifies parts of the physics
(coarser description than in \cite{Zhou-2001-ERM}), but includes the essence
of atmosphere-fire interaction.
%end copied
Future developments of the numerics and parallelization of our model are
expected to be able to handle realistic fires of the size of several $km$,
coupling with the atmosphere, and assimilation of real-time data.

An important point is that our paradigm attempts to strike a balance between
model complexity and fast execution. Thus, the present model is based on just
two PDEs in two horizontal spatial dimensions: prognostic (predictive)
equations for heat and fuel. We use a single semi-empirical reaction rate to
achieve the desired combustion model results. In other words, we solve the set
of equations known as the reaction-convection-diffusion problem using reaction
rates based on the Arrhenius equation, which relates the rate at which a
chemical reaction proceeds to the temperature.

This is the simplest combustion model and it is known to produce solutions
with traveling combustion waves, that is, a propagating area of localized
combustion made up of the preheated area ahead of the fire, the combustion
zone, and the post-frontal burning region. One reason for considering a
PDE-based model is that even simple reaction-diffusion equations are capable
of the complex nonlinear, unsteady behavior such as pulsation and bifurcation
that is seen in reality but cannot be reproduced by empirical models.

The characteristics of the combustion wave (maximal temperature, width of the
burning region as defined by the leading and trailing edge, and the speed of
propagation)\ are used to calibrate the parameters of the model. We note that
physical behavior can be achieved by a very simple model that can reproduce
realistic fire behavior very quickly on today's computers. The PDE-based
models in this paper are not necessarily original, cf., e.g.,
\cite{Weber-1997-CWG}, and some PDE coefficients have been determined empirically
from measured time-temperature curves before
\cite{Balbi-1999-DMF,Giroud-1998-DAT,Morandini-2001-CRH}, though not in a
reaction-diffusion PDE model like the model here. We provide a new systematic
procedure for the calibration of the PDE model on real wildland fire data
based on separation of nondimensional properties and solution scales. The
general calibration problem is an interesting optimal control and stochastic
parameter estimation problem in its own right that will be studied in detail
in future works.

We then proceed to data assimilation to modify the state of the running model
from data. A version of the ensemble Kalman filter (EnKF) is used. This work
appears to be the first wildland fire model with data assimilation.

Future extensions of this work include coupling with a numerical weather
prediction model, modeling of water content in the fuel, multiple fuel layers,
separate treatment of the gas phase (i.e. pyrolysis), crown fire, modeling
spotting by stochastic differential equations, preheating by long-range
radiation, and contemporary numerical methods such as finite elements and
level set methods. The exact number of data points (in space and time) that
are necessary for recovering good predictions when the numerical solution
diverges from reality depends on the particular data assimilation method used.
For the method in this paper, the correction in the location of the fireline
should not be larger than the width of the reaction zone. This is similar to
the situation for data assimilation into hurricane models, where the
correction in the location of the vortex should not be larger than the vortex
size \cite{Chen-2007-AVP}. Advanced data assimilation methods that allow
sparser data and larger correction are the subject of further research. While
real-time data are routinely available for weather forecasting systems, in
a~wildland fire the data collection is less straightforward. Available data
include multi-spectrum infrared airborne photographs, processed to recover the
fire region and to some extent the temperature, and radioed data streams from
hardened sensors put in the fire path
\cite{Kremens-2003-MTT,Ononye-2005-IFT,Ononye-2007-AEF}. For overviews of the
whole project including computer science aspects, data collection, and
visualization, see \cite{Mandel-2004-NDD,Mandel-2005-TDD,Mandel-2007-DDD} and
\cite{Douglas-2006-DVW}.

The remainder of the paper is organized as follows. In
Section~\ref{sec:formulate-model}, we state our PDE-based fire model. Then in
Section~\ref{sec:related-models}, we describe the relation of our model to
other models in the literature. In Section~\ref{sec:derive-model}, the model
is derived from physical principles in more detail. In
Section~\ref{sec:identification}, we develop a method to determine the
coefficients of the PDE model using wildland fire data. In
Section~\ref{sec:enkf}, we describe the ensemble Kalman filter techniques for
the data assimilation. In Section~\ref{sec:numerical-results}, we test the PDE
and ensemble Kalman filter methods on a two-dimensional representation of a
wildland fire and calibrate the models against real data. Finally,
Section~\ref{sec:conclusions} contains our conclusions.

\section{Formulation of the model}

\label{sec:formulate-model} We consider the model of fire in a layer just
above the ground.
First, we define the following terms:
\begin{description}
\item $T$ ($K$) is the temperature of the fire layer,

\item $S\in\left[  0,1\right]  $ is the fuel supply mass fraction (the
relative amount of fuel remaining),

\item $k$ ($m^{2}s^{-1}$) is the thermal diffusivity,

\item $A$ ($Ks^{-1}$) is the temperature rise per second at the maximum
burning rate with full initial fuel load and no cooling present,

\item $B$ ($K$) is the proportionality coefficient in the modified Arrhenius law,

\item $C$ ($K^{-1}$) is the scaled coefficient of the heat transfer to the
environment, $\ $

\item $C_{S}$ ($s^{-1}$) is the fuel relative disappearance rate,

\item $T_{a}$ ($K$) is the ambient temperature, and

\item $\overrightarrow{v}$ ($ms^{-1}$) is the wind speed given by atmospheric
data or model.
\end{description}
The model is derived from the
conservation of energy, balance of fuel supply, and the fuel reaction rate:%

\begin{align}
\frac{dT}{dt} &  =\nabla\cdot\left(  k\nabla T\right)  -\overrightarrow
{v}\cdot\nabla T+A\left(  Se^{-B/(T-T_{a})}-C\left(  T-T_{a}\right)  \right)
,\label{eq:heat}\\
\frac{dS}{dt} &  =-C_{S}Se^{-B/(T-T_{a})},\quad T>T_{a},\label{eq:fuel}%
\end{align}
with the initial values%
\begin{equation}
S\left(  t_{\text{init}}\right)  =1\text{ and }T\left(  t_{\text{init}%
}\right)  =T_{\text{init}}\text{.}\label{eq:initial}%
\end{equation}

The diffusion term $\nabla\cdot\left(  k\nabla T\right)  $ models short-range
heat transfer by radiation in a semi-permeable medium, $\overrightarrow
{v}\cdot\nabla T$ models heat advected by the wind, $Se^{-B/(T-T_{0})}$ is the
rate fuel is consumed due to burning, and $AC\left(  T-T_{a}\right)  $ models
the convective heat lost to the atmosphere. The reaction rate $e^{-B/(T-T_{a}%
)}$ is obtained by modifying the reaction rate $e^{-B/T}$ from the Arrhenius
law by an offset to force zero reaction at ambient temperature, with the
resulting reaction rate smoothly dependent on temperature.

A more detailed derivation of the model from physical principles is contained
in Section~\ref{sec:derive-model}. Calibration of the coefficients from
physically observable quantities will be described in
Section~\ref{sec:identification}.

\section{Relation to other models}

\label{sec:related-models}

Recent surveys of wildland fire models and their histories are in
\cite{Morvan-2002-BWF,Pastor-2003-MMC,Seron-2005-EWF}.

\subsection{Models based on diffusion-reaction PDEs}

It is known that systems of the form (\ref{eq:heat}-\ref{eq:fuel}) admit
traveling wave solutions. The temperature in the traveling wave has a sharp
leading edge, followed by an exponentially decaying cool-down trailing edge.
This was observed numerically but we were not able to find a rigorous proof in
the literature in exactly this case, though proofs for some related systems exist.

For a related system with fuel diffusion, the existence and speed of traveling
waves were obtained by asymptotic methods already in classical work,
summarized in the monograph \cite{Zeldovich-1985-MTC}. For the system
(\ref{eq:heat}-\ref{eq:fuel}), \cite{Weber-1997-CWG} obtains approximate
combustion wave speed by heuristic asymptotic methods under the assumption
that no heat is lost and ambient temperature is absolute zero, which is
equivalent to our setting $C=0$. Models that do not guarantee zero combustion
at ambient temperature suffer from the \textquotedblleft{}cold boundary
difficulty\textquotedblright: by the time a combustion wave gets to a given
location, the fuel at that location is depleted by the ongoing reaction at
ambient temperature. So, no perpetual traveling combustion waves can exist,
and there are only \textquotedblleft pseudo-waves\textquotedblright\ that
travel only for a finite time
\cite{Berestycki-1991-MIC,Mercer-1996-CPW,Zeldovich-1985-MTC}.
\cite{Campos-2004-RDP} derives the speed of traveling waves in a
simplified model with the reaction started by ignition at a~given temperature,
followed by an assumed temperature indepentent reaction rate that is only
proportional to the fuel remaining. This is similar to the model in
\cite{Balbi-1999-DMF}.

Equation (\ref{eq:heat}) without fuel depletion (i.e., with constant $S$) and
without wind (i.e., $\overrightarrow{v}=0$) is a special case of the nonlinear
reaction-diffusion equation,%
\begin{equation}
\frac{dT}{dt}=\nabla\cdot\left(  \nabla T\right)  +f\left(  T\right)  .
\label{eq:reaction-diffusion}%
\end{equation}
Reaction-diffusion equations of the form of equation
(\ref{eq:reaction-diffusion}) are known to possess traveling wave solutions,
which switch between values close to stationary states given by $f\left(
T\right)  =0$ \cite{Gilding-2004-TWN,Infeld-2000-NWS}. The simplest model
problem is Fisher's equation with $f\left(  T\right)  =T\left(  1-T\right)  $,
for which the existence of a traveling wave solution and a~formula for its
speed were proven in \cite{Kolmogorov-1937-SDE}. For an analytical study of
the evolution of waves to a traveling waveform, see \cite{Sherratt-1998-TID},
and for a numerical study, see \cite{Gazdag-1974-NSF,Zhao-CDS-2003}%
. \cite[Ch. 8 and 11]{Robinson-2001-IDP} gives proofs of the existence of
solution and attractors for reaction-diffusion equations, but does not mention
traveling waves. \cite{Asensio-2002-WFM} proves the existence of a solution,
but not traveling waves, for the reaction-diffusion equation (\ref{eq:heat})
except with a nonlinear diffusion term $\nabla\cdot\left(  kT^{3}\nabla
T\right)  $ and again without considering fuel depletion.

\cite{Seron-2005-EWF} considers a two reaction model (solid and pyrolysis gas),
and argues that the modeling of pyrolysis by a separate reaction is essential
for capturing realistic fire behavior. For a more complicated model of this
type that includes other components, e.g., water vapor, see \cite{Grishin-1996-GMM}.

Various aspects of special cases of equations (\ref{eq:heat}-\ref{eq:fuel})
have been studied in a~number of papers. \cite{Weber-1991-TCW} uses a formal
expansion in an Arrhenius reaction model to get the wave speed and a
prediction of whether a small fire will or will not spread.
\cite{Mercer-1995-CWS} computes the ignition wave speed and extinction wave
speed numerically. The speed and stability of combustion waves are analyzed
by asymptotic expansion in \cite{Gubernov-2003-EFS}. An approximation to the
temperature reaction equation gives the size of the reaction zone and the slope of
the temperature curve. \cite{Norbury-1988-TCP,Norbury-1988-TCW} derive a
nonlinear eigenvalue problem for a traveling wave in a different combustion
problem, with fuel reaction, solve it numerically by the shooting method, and
study the existence and stability of the traveling wave solution. See also
\cite{Gubernov-2003-EFS,Gubernov-2004-EFS} for the gas case (i.e., also with
fuel diffusion instead of just temperature diffusion).

Enriched finite element methods for the linear diffusion-advection-reaction
problem are designed and error estimates given in
\cite{Asensio-2004-RBM,Codina-1998-CFE,Franca-2005-TMS,Franca-2006-EFM}.
However, in all those works, the reaction function $f$ is replaced by a linear
function, so there are no traveling wave solutions, and fuel consumption is
not considered. \cite{Mickens-2005-NFD} compares several time discretization
techniques in the presence of nonlinear reaction terms.
\cite{Asensio-2000-TEE,Ferragut-2002-MFE} provide error estimates
for mixed finite elements applied to the single species combustion equation
and more general reaction-diffusion problems, but again without a fuel balance
equation. \cite{Sembera-2001-NGM} proposes
a nonlinear Galerkin method for reaction-diffusion problems, and proves
convergence by a compactness argument. Even simple nonlinear models exhibit
bifurcations, which can be examined by direct simulation
\cite{Theodoropoulos-2000-CSB}. Approximation of Fisher's equation by finite
elements are studied numerically in \cite{Carey-1995-LFE,Roessler-1997-NSD},
especially regarding the correct wave speed, but no
error estimates are given. For finite differences applied to a
reaction-diffusion equation, see \cite{Liao-2006-FCA}. Since the common
feature of the solutions of reaction-diffusion equations is the development of
a sharp wave, with the solution being almost constant elsewhere, interface
tracking techniques such as level set methods
\cite{Sethian-1999-LSM,Sussman-1994-LSA} are relevant here as well.

\subsection{Fireline evolution, fire spread, and empirical models}

The reaction zone in reaction-diffusion models is typically very thin, and
resolving it correctly requires very fine meshes. Hence, a number of models
consider the evolution of the fireline instead. Combustion equations in the
reaction sheet limit or large activation energy asymptotics reduce
to a representation of the reaction zone (here, the fireline) as an evolving
internal interface
\cite{Chen-1992-GPIs,Class-2003-UMF,Dold-2003-HOE,Fife-1988-DIL,Law-1993-CAE,Rastigejev-2006-NSF}%
, though this reduction does not seem to have been done for exactly the same
equations as here. The asymptotic models typically compute the speed of the
movement of the reaction interface in the normal direction, often involving
its curvature.

Fireline evolution models often postulate empirically observed properties of
the fire, such as the fire spread rate in the normal direction, instead of
physically based differential equations. Modeling of fireline evolution was
reviewed by \cite{Weber-1991-MFS}.
\cite{Albini-1994-PBS,Albini-1997-ICL,Frandsen-1971-SFP}
derive fire spread rates
without using reaction kinetics. \cite{Richards-1995-GMF,Richards-1999-MMC}
study evolution of the fireline as a curve. \cite{Dupuy-1999-FSP} introduces a convective term in a
radiation based model in an attempt to better describe slope effects on rate
of spread.

Far fewer examples exist where data was used to calibrate models.
Fire spread models based on radiation like \cite{Albini-1986-WFS}
are tested against data in \cite{Dupuy-2000-TTR}.
\cite{Albini-1994-PBS,Albini-1997-ICL} calibrate a mass loss model
by crib burning experiments. \cite{Rothermel-1972-MMP} formulate a model for
surface fire spread rate with a physically based core rate of spread in zero
wind on flat ground, calibrated to other wind speed and slopes using
laboratory measurements. \cite{Balbi-1999-DMF} use the energy balance
equation coupled with a model of ignition at a threshold temperature, followed
by exponential decay. Measured temperature profiles are used to identify
parameters of the model. \cite{Morandini-2001-CRH,Simeoni-2001-WIF}
use laboratory data to validate predictions made with different systems of
PDEs. \cite{Wotton-1999-EFF} use actual fire spread data and theoretical
calculations to test the effect of fire front width on surface spread rates
through radiative transfer terms. \cite{Viegas-2005-MMF} postulates empirical
rates of fire spread and of the wind created by the fire and identifies the
coefficients from experiments. The feedback between the fire and the
surrounding flow is then modeled by a simple one-dimensional differential
equation, which is sufficient to explain the conditions for the fire spread to
stop or accelerate to a blowup.

\subsection{Coupled fluid-fire models}

Wildland fire models (either empirical, semi-empirical, or PDE-based) have
been coupled to a fluid environment that may be (for small domains) a
computational fluid dynamics model (i.e., models the flow and thermodynamics of
air) or a numerical weather prediction model (i.e., a computational fluid
dynamics model that also considers moist atmospheric processes, the formation
of precipitation, and flow over topography).

The FIRETEC model \cite{Linn-2002-SWB,Linn-1996-FTM} simulates wildland fires
by representing the average reaction rates and transport over a resolved
volume, usually on the order of $1 m^{3}$ in three dimensional space. This
attempts to resolve the effects of heat transfer processes without
representing each in detail. The ambient environment is air, but the model
omits weather processes. \cite{Larini-1998-MFF} gives a multiphase, reactive,
and radiative one dimensional model specialized for wildland fires.
\cite{Dupuy-2005-NSC} adds detailed fluid modeling to the model from
\cite{Larini-1998-MFF} to study crown fires in air, again with no weather
processes and in two dimensions. \cite{Grishin-2002-MMS} presents a complicated
model of surface and canopy fire based on fluid dynamics and chemical reaction
equations where prognostic equations are created for gases that have been
grouped into reactive combustible gases, non-reactive combustion products, and
an oxidizer ($O_{2}$). Methods of this type are standard in combustion
modeling. For a flame model with detailed chemistry and physics, see
\cite{Ern-1995-DCM}.

An alternate approach is adopted in \cite{Clark-2004-DCA,Clark-1996-CAF},
where a semi-empirical fire spread model based upon the
\cite{Rothermel-1972-MMP} fire spread equation and a canopy fire model are
coupled to a numerical weather prediction model to model the interactions
between wildland fires and the atmospheric environment. Here, weather
processes ranging from synoptic to boundary layer scale are simulated with
good fidelity, and the combustion processes are represented by
semi-empirical formulas in order to capture the sensible (temperature) and
latent (water vapor) heat fluxes into the environment.

\section{Derivation of the model}

\label{sec:derive-model}

We consider fire in a ground layer of some unspecified finite small thickness
$h$. The fire layer consists of the fuel and air just above the fuel. All
modeled quantities are treated as two dimensional, homogenized in the vertical
direction over the ground layer. We will not attempt to derive equations and
substitute coefficients from material properties because of the degree of
simplification and uncertainty present in the homogenization. Instead,
physical laws will be used to derive the form of the equations and the
coefficients will be identified later from the dynamical behavior of the
solution. We first derive the system of PDEs based on conservation of energy
and fuel reaction rate in Section \ref{sec:conservation} and then discuss the
choice of the reaction term in Section \ref{sec:reaction}.

\subsection{Heat and fuel supply balance equations}

\label{sec:conservation}

The chemical reactions are a heat source. Heat transfer is due to radiation
and convection to the atmosphere. The short-range heat transfer due to
radiation and turbulence is modeled by diffusion. The two dimensional heat
flux through a segment per length unit then is
\begin{equation}
\overrightarrow{q}_{r}=-k_{1}\nabla T\quad(Wm^{-1}). \label{eq:radiation-flux}%
\end{equation}
The constant $k_{1}$ ($WK^{-1}$) will be identified later.

Heat is generated by the chemical reaction of burning. We model the burning as
a reaction in which the rate depends on temperature only, so the reaction rate
is $C_{S}r\left(  T\right)  $, where $C_{S}$ is a coefficient of
proportionality ($1/s$), and $r$ is dimensionless. Let $F>0$ ($kg\,m^{-2}$) be
the concentration of fuel remaining. Then the rate at which the fuel is lost
is proportional to the rate of reaction and the amount of fuel available,%
\begin{equation}
\frac{dF}{dt}=-FC_{S}r\left(  T\right)  . \label{eq:fuel-lost}%
\end{equation}
The heat generated per unit surface area is then proportional to the fuel
lost,
%\begin{align}
%\label{e4.7}q_{g} = A_{1} S e^{-B/(T-T_{0})}%
%\end{align}
%{\bf Propose change to%
\begin{equation}
q_{g}=A_{1}FC_{S}r\left(  T\right)  ,\quad(Wm^{-2}) \label{eq:heat-generated}%
\end{equation}
where $A_{1}$ ($J\;kg^{-1}$) is the heat released per unit mass of fuel.

Heat per unit area lost due to natural convection to the atmosphere due to
buoyancy is given by Newton's law of cooling,%
\begin{equation}
q_{c}=C_{a}(T-T_{a}),\quad\left(  Wm^{-2}\right)  \label{eq:heat-lost}%
\end{equation}
where $T_{a}$ is the ambient temperature ($K$) and $C_{a}$ ($Wm^{-2}K^{-1}$)
is the heat transfer coefficient. In this model, it is assumed that the
convective heat transfer is dominant, and so the effect of radiation into the
atmosphere is included in (\ref{eq:heat-lost}).

The material time derivative of the temperature is given by
\begin{equation}
\frac{DT}{Dt}=\frac{dT}{dt}+\overrightarrow{v}\cdot\nabla T .
\label{eq:material}%
\end{equation}
$dT/dt$ is the Eulerian
(or spatial) time derivative of temperature, $\overrightarrow{v}$ is the
(homogenized) velocity of the air, $\rho$ is the homogenized surface density
of the fire layer ($kg\;m^{-3}$), and $c_{p}$ is the homogenized specific heat
of the fire layer ($J\,kg^{-1}K^{-1}$). Again, none of the coefficients $h$,
$\rho$, or $c_{p}$ can be assumed to be actually known. The units of the
product $h\rho c_{p}$ are $JK^{-1}m^{-2}$.

From the divergence theorem, we now obtain the conservation of energy in the
fire layer as%
\begin{equation}
h\rho c_{p}\frac{DT}{Dt}=\nabla\cdot\overrightarrow{q}_{r}+q_{g}-q_{c}.
\label{eq:heat-conservation}%
\end{equation}

The velocity vector $\overrightarrow{v}$ is obtained from the state of the
atmosphere as data, or in future work by coupling with an atmospheric model.
In the present model, the velocity vector also incorporates the effect of
slope: adding to the wind a small multiple of the surface gradient somewhat
simulates the effect that fire spreads more readily uphill. In addition, since
the speed of the air is zero at the ground level if as usual no-slip
conditions are assumed, the homogenized speed through the fire layer should be
approximated by scaling the given wind velocity by a constant less than one.

We now write the equations in a form suitable for identification of the
coefficients (which will be formally done in
Sec. \ref{sec:identification}). The goal is to obtain a system
of equations with a minimal number of coefficients in as simple form as
possible. In addition, we wish to relate the coefficients to the behavior of
the solution for certain particular coefficient values, rather than to
material and physical properties of the medium, which are in general unknown.
We introduce the mass fraction of fuel by%
\[
S=\frac{F}{F_{0}},
\]
where $F_{0}$ is the initial fuel quantity. Substituting in the appropriate
values for the heat sources and fluxes, (\ref{eq:radiation-flux}),
(\ref{eq:heat-generated}), and (\ref{eq:heat-lost}), into (\ref{eq:material})
and (\ref{eq:heat-conservation}), and some simple algebra, we obtain the
energy balance and fuel reaction rate equations%

\begin{align}
\frac{dT}{dt} &  =\nabla\cdot\left(  k\nabla T\right)  -\overrightarrow
{v}\cdot\nabla T+A\left(  Sr\left(  T\right)  -C_{0}\left(  T-T_{a}\right)
\right)  ,\label{eq:heat-r}\\
\frac{dS}{dt} &  =-C_{S}Sr\left(  T\right)  ,\label{eq:fuel-r}%
\end{align}
with%
\[
k=k_{1}/(h\rho c_{p}),\quad A=A_{1}C_{S}/(hc_{p}\rho),\quad\mbox{and}%
\quad C_{0}=C_{a}/A_{1}.
\]

Alternatively, we could have taken disappearance of fuel on the left hand side
of the heat balance equation (\ref{eq:heat-conservation}) (fuel that has
burned does not need to be heated), which would lead to an equation of the
form%
\[
\left(  1+C_{1}S\right)  \frac{dT}{dt}=\nabla\cdot\left(  k\nabla T\right)
+\overrightarrow{v}\cdot\nabla T+A\left(  Sr\left(  T\right)  -C_{0}\left(
T-T_{a}\right)  \right)
\]
instead of (\ref{eq:heat-r}). We have chosen not to do so since our goal is
to work with the simplest possible model, whose coefficients can be
identified. The fuel disappearance in the heat equation affects the
temperature profile significantly only in the reaction zone, which is the
highest part of the temperature curve. Before the ignition, there is a full
fuel load, $S=1$, and after a fairly short reaction time, well before most of
the cooling takes place, the remaining fuel settles to some residual value
which then remains constant. The effect of the decreased heat capacity of the
remaining fuel is then absorbed into the cooling term $AC_{0}$.

\subsection{Reaction rate}

\label{sec:reaction}

The Arrhenius reaction rate from physical chemistry is given by
\begin{equation}
r\left(  T\right)  =e^{-B/T}, \label{eq:arrhenius}%
\end{equation}
where the coefficient $B$ has units $K$. This equation is valid only for gas
fuel premixed with a sufficient supply of oxygen. This approximation ignores
fuel surface effects but it is widely used nonetheless. One consequence of
(\ref{eq:arrhenius}) is that the reaction has a nonzero rate at any
temperature above absolute zero. Since the time scale for burning is much
smaller than the oxidation rate at ambient temperature, we modify
(\ref{eq:arrhenius}) so that no oxidation occurs below some fixed temperature,
$T_{0}$ (see Section \ref{sec:identification}), and take instead%
\begin{equation}
r\left(  T\right)  =\left\{
\begin{array}
[c]{c}%
e^{-B/\left(  T-T_{0}\right)  },\quad T>T_{0},\\
0,\quad T\leq T_{0}.
\end{array}
\right.  \label{eq:rate}%
\end{equation}
Note that the fuel consumption rate is a smooth function of $T$, which is
favorable for a numerical solution, unlike in \cite{Asensio-2002-WFM}, where a
cutoff function was used.

\section{Identification of coefficients}

\label{sec:identification}

We wish to use an observed behavior of the fire rather than physical material
properties to identify the coefficients. It is not simple to obtain reasonable
behavior of the solution from substituting physical coefficients into the
equations. Further, as explained in Section \ref{sec:derive-model}, because of
a number of simplifying assumptions employed and because of the homogenization
of coefficients over a fire layer of unspecified thickness, it is not quite
clear what the material properties should be anyway.

We first use basic reaction dynamics and a reduced model to find rough
approximate values of the coefficients that produce a reasonable solution.
Then we transform the equations to a nondimensional form, which allows us to
separate the coefficients into those that determine the qualitative behavior
of the solution and those that determine the scales. We propose to use the
approximate coefficients obtained from the reduced models as initial values
for identification of the coefficients by the nondimensionalization method to
match observed temperature profiles.

\subsection{Reaction rate coefficients}

Consider first the hypothetical case in which $T$ is
constant in space, so that only heat due to the reaction (burning) and natural
convection contribute non-zero terms in the heat equation, (\ref{eq:heat}),
and essentially the full initial fuel supply $F_{0}$ is present at all times
(the rate of fuel consumption is negligible, $C_{S}\approx0$, so $S\approx
1$):
\begin{equation}
\frac{dT}{dt}=A\left(  e^{-B/\left(  T-T_{0}\right)  }-C\left(  T-T_{a}%
\right)  \right)  . \label{eq:no-fuel-loss}%
\end{equation}
Constant values of temperature which are solutions to (\ref{eq:no-fuel-loss})
are called equilibrium points, and at these points the heat produced by the
reaction equals the heat lost to the environment,%
\begin{equation}
f(T)=e^{-B/\left(  T-T_{0}\right)  }-C\left(  T-T_{a}\right)  =0.
\label{eq:balance}%
\end{equation}

\begin{figure}[ptb]
\begin{center}
\includegraphics[width=3in]{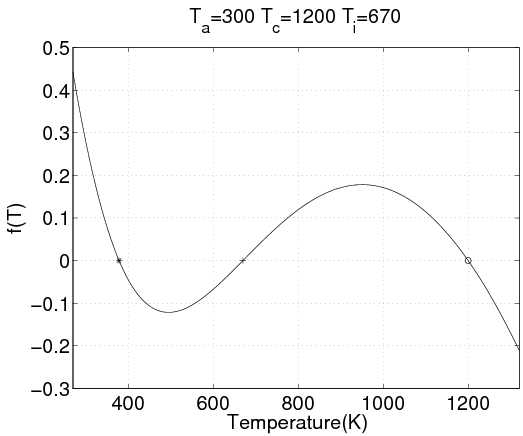}
\end{center}
\par
.\caption{Sample reaction heat balance function $f(T)$ from equation
(\ref{eq:balance})}%
\label{fig:f}%
\end{figure}

\begin{figure}[ptb]
\begin{center}
\includegraphics[width=3in]{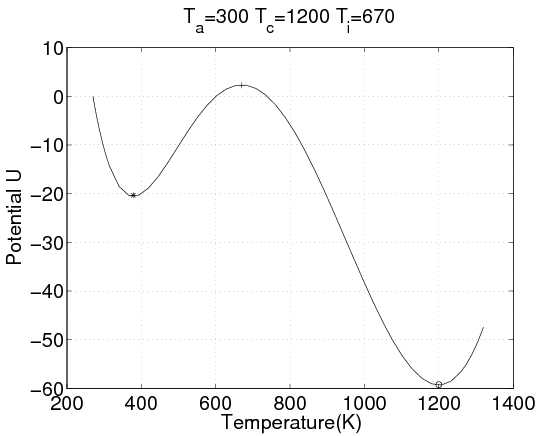}
\end{center}
\caption{Reaction heat balance potential $U$}%
\label{fig:U}%
\end{figure}

Equation (\ref{eq:balance}) has at most three roots
\cite{Frank-Kamenetskii-1955-DHE}, see for example,~Figure \ref{fig:f}. The
first zero, denoted as $T_{p}$, called the lower temperature regime by
\cite{Frank-Kamenetskii-1955-DHE}, is a stable equilibrium temperature. If
the temperature goes below this temperature then the heat generated from the
reaction dominates and the temperature rises. If the temperature goes above
this temperature then convective cooling dominates and the temperature
decreases. If $T_{0}<T_{a}$, then this point is typically just above the
ambient temperature, $T_{a}$, since some reaction is present even at ambient
temperature. The middle zero, $T_{i}$, is an unstable equilibrium point. If
the temperature goes below $T_{i}$ then convective cooling dominates and the
temperature decreases. Above $T_{i}$, the heat due to chemical reactions
dominate and the temperature increases. We refer to $T_{i}$ as the
\emph{auto-ignition temperature}, the temperature above which the reaction is
self-sustaining \cite{Quintiere-1998-PFB}. The stable equilibrium at a high
temperature, $T_{c}$, is the maximum stable combustion temperature, assuming
replenishing of the supply of fuel and oxygen. The temperature $T_{c}$ is
called the high temperature regime by \cite{Frank-Kamenetskii-1955-DHE}. The
stability properties of the equilibrium points are also clear from the graph
of the potential $U(T),$ defined by $U^{\prime}(T)=f(T)$. The stable
equilibrium points are local minima of the potential, while the autoignition
temperature is a local maximum and thus an unstable equilibrium
(see Fig.~\ref{fig:U}).

While the coefficients $B$ and $C/F_{0}$ in (\ref{eq:balance}) are generally
unknown, they can be found from the equilibrium temperature points. Suppose
that two roots $T_{i}$ and $T_{c}$ of $f(T)$ from (\ref{eq:balance}) are given
such that $T_{0}\leq T_{a}<T_{i}<T_{c}$. Then simple algebra gives
\begin{equation}
B=\frac{\ln{\left(  \frac{T_{i}-T_{a}}{T_{c}-T_{a}}\right)  }}{\frac{1}%
{T_{c}-T_{0}}-\frac{1}{T_{i}-T_{0}}}\quad\mbox{and}\quad%
C=\frac{e^{-B/(T_{i}-T_{0})}}{T_{i}-T_{a}}. \label{eq:BC}%
\end{equation}

It should be noted that the coefficients $B$ and $C$ from (\ref{eq:BC}) result
in three equilibrium points only when $T_{i}$ is significantly higher than
$T_{a}$. When $T_{i}$ is too close to $T_{a}$, the resulting energy balance
equation $f(T)=0$ has only two roots. This, however, does not occur for the
values of $T_{a}$, $T_{i}$, and $T_{c}$ of interest.

\begin{figure}[ptb]
\begin{center}
\includegraphics[width=3in]{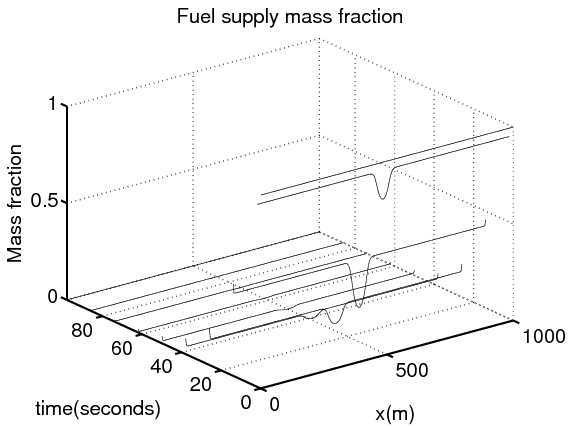} \includegraphics[width=3in]{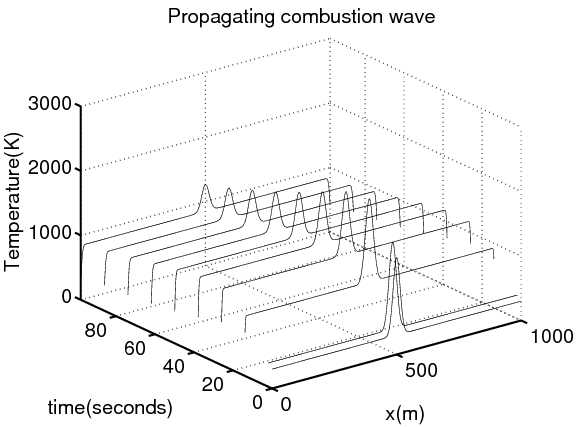}
\end{center}
\caption{Solution with the Arrhenius reaction rate. Due to nonzero reaction
rate at ambient temperature, fuel starts disappearing and thus a propagating
combustion wave does not develop. The coefficients are $k=2.1360\times
10^{-1}m^{2}s^{-1}K^{-3}$, $A=1.8793\times10^{2}Ks^{-1}$, $B=5.5849\times
10^{2}K$, $C=4.8372\times10^{-5}K^{-1}$, and $C_{S}=1.6250\times10^{-1}s^{-1}%
$. $T_{c}=1200K$ and $T_{i}=670K$.}%
\label{fig:f0}%
\end{figure}

\begin{figure}[ptb]
\begin{center}
\includegraphics[width=3in]{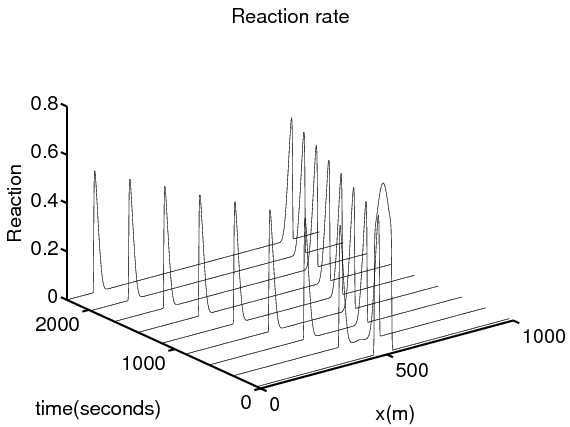} \includegraphics[width=3in]{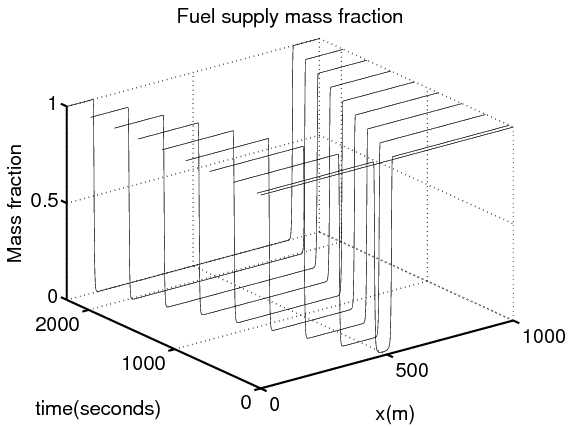}
\includegraphics[width=3in]{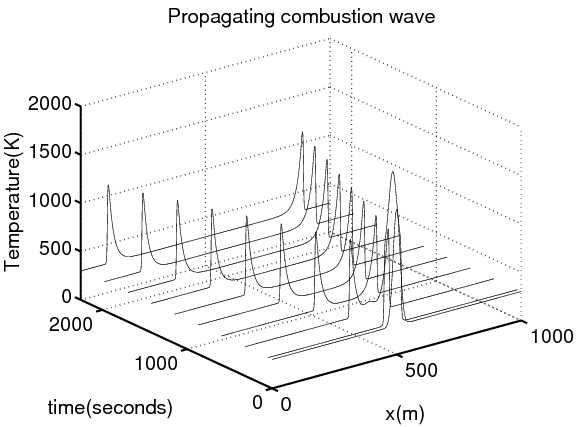}
\end{center}
\caption{Solution with Arrhenius reaction rate modified by temperature offset
$T_{a}$ to force zero reaction rate at ambient temperature. A propagating
combustion wave develops. The coefficients are $k=2.1360\times10^{-1}%
m^{2}s^{-1}K^{-3}$, $A=1.8793\times10^{2}Ks^{-1}$, $B=5.5849\times10^{2}K$,
$C=4.8372\times10^{-5}K^{-1}$, and $C_{S}=1.6250\times10^{-1}s^{-1}$.
$T_{c}=1200K$ and $T_{i}=670K$.}%
\label{fig:ra}%
\end{figure}

First consider the solution of (\ref{eq:heat-r}-\ref{eq:fuel-r}) and the
reaction rate (\ref{eq:rate}) with $T_{0}=0$. Then the reaction is the
Arrhenius rate known from chemistry and there is a nonzero reaction rate at
$T=T_{a}$. This results in fuel loss everywhere, and, in our computational
experiments, no traveling combustion wave developed (see Fig.~\ref{fig:f0}) since,
after a relatively short time, there was not enough fuel to sustain
combustion. This phenomenon is known as the cold boundary effect in combustion
literature \cite{Weber-1997-CWG}: a traveling combustion wave solution does
not exist, and there can only be pseudo-waves that propagate for a limited
time \cite{Berestycki-1991-MIC,Mercer-1996-CPW} and then vanish. Since for the
values of $B$ and $C$ obtained from realistic $T_{i}$ and $T_{c}$, the fuel
disappears rather quickly, we force the reaction rate, $r(T)$, to be zero at
ambient temperature by choosing the offset $T_{0}=T_{a}$. Using the offset by
$T_{a}$ is essentially the same as assuming that the ambient temperature is
absolute zero as commonly done in combustion literature \cite{Weber-1997-CWG}.
In this case, a propagating combustion wave develops (see Figs.~\ref{fig:ra}
and \ref{fig:tw}). Therefore, we use $T_{0}=T_{a}$.

It should be noted that traveling combustion waves, such as in Fig.
\ref{fig:tw}, are caused by the combined effect of reaction and diffusion;
convection does not play a role in this section. The reaction heat diffuses
forward on the leading edge, heating the fuel ahead of the wave, until the
reaction ahead of the wave can sustain itself, thus causing the combustion to
spread. On the trailing edge, the reaction drops off due to fuel depletion,
after which temperature decays due to cooling.

\subsection{Cooling coefficient}

\label{sec:ABC}

The coefficients $B$ and $C$ in the modified reaction form (\ref{eq:balance})
 have been determined from reasonable values of $T_{c}$ and $T_{i}$ by
(\ref{eq:BC}). We want to determine the remaining coefficient, $A$. This can
be done from the characteristic cooling time. Consider the trailing edge of a
traveling combustion wave, after all or most of the fuel has been depleted,
temperature drops, and heat generated by the reaction and diffusion drop to an
insignificant level. From that point on, the temperature
approximately satisfies%
\[
\frac{dT}{dt}=-AC(T-T_{a}).
\]
Thus, at the trailing edge, given $T$ at some time $t_{0}$, we have
\[
T(t)=T_{a}+\left(  T(t_{0})-T_{a}\right)  e^{-AC(t-t_{0})}%
\]
and we can define the characteristic cooling time, $t_{c}$, to be the time
which the fire layer takes to cool by a factor of $1/e$, i.e.,%
\[
T(t_{0}+t_{c})-T_{a}=\frac{1}{e}\left(  T(t_{0})-T_{a}\right)  ,
\]
which proves that$ACt_{c}=1$, or%
\begin{equation}
A=\frac{1}{Ct_{c}}. \label{eq:A}%
\end{equation}

\subsection{Scales and non-dimensional coefficients}

\label{sec:non-dimensional}

We now write the model in terms of nondimensional variables, which control the
qualitative behavior of the system (\ref{eq:heat}-\ref{eq:fuel}). Again, we do
not consider the wind here yet, and so
\begin{align}
\frac{dT}{dt}  &  =\nabla\left(  k\nabla T\right)  +A\left(  Se^{-\frac
{B}{T-T_{a}}}-C\left(  T-T_{a}\right)  \right)  ,\label{eq:heat2}\\
\frac{dS}{dt}  &  =-SC_{S}e^{-\frac{B}{T-T_{a}}},\qquad T>T_{a}.
\label{eq:fuel2}%
\end{align}

The substitution%
\[
\widetilde{T}=\frac{T-T_{a}}{B},\qquad\widetilde{x}=\frac{x}{k^{1/2}%
B^{1/2}A^{-1/2}},\qquad\mbox{and}
\qquad\widetilde{t}=\frac{tA}{B}
\]
transforms (\ref{eq:heat2}-\ref{eq:fuel2}) into a non-dimensional form
\begin{align}
\frac{d\widetilde{T}}{d\widetilde{t}}  &  =\widetilde{\nabla}\cdot\left(
\widetilde{\nabla}\widetilde{T}\right)  +\widetilde{S}e^{-1/\widetilde{T}%
}-\lambda\widetilde{T},\label{eq:dimensionless-heat}\\
\frac{d\widetilde{S}}{d\widetilde{t}}  &  =-\beta\widetilde{S}e^{-1/\widetilde
{T}},\qquad\widetilde{T}>0, \label{eq:dimensionless-fuel}%
\end{align}
with two dimensionless coefficients%
\begin{equation}
\lambda=CB\quad\mbox{and}\quad%
\beta=\frac{BC_{S}}{A}. \label{eq:dimensionless-parameters}%
\end{equation}

Therefore, the qualitative behavior of the solution is determined only by the
nondimensional coefficients $\lambda$ and $\beta$, which can be varied independently.

The nondimensional form (\ref{eq:dimensionless-heat}%
-\ref{eq:dimensionless-fuel}) suggests a strategy for identification of the
coefficients $k$, $A$, $B$, $C$, $C_{S}$:\ first match nondimensional
properties of the traveling combustion wave, such as the ratio of the width of
the leading edge and the trailing edge and the fuel fraction remaining after
the combustion wave by varying $\lambda$ and $\beta$. The width of the wave
can be measured e.g. as the distance of the points where the temperature
equals 50\% of the maximum. The nondimensional traveling wave solution
$\widetilde{T}(\widetilde{t},\widetilde{x})$, $\widetilde{S}(\widetilde
{t},\widetilde{x})$ has some (nondimensional) maximal temperature
$\widetilde{T}_{\max}$, width $\widetilde{w}$, and speed $\widetilde{v}$,
while the data (a measured temperature profile) has maximal temperature
$T_{\max}$, the width $w$, and the speed $v$ of the traveling wave. This
determines the scales%
\[
T_{1}=\frac{T_{\max}}{\widetilde{T}_{\max}},\quad x_{1}=\frac{w}{\widetilde
{w}},\quad\mbox{and}\quad%
t_{1}=\frac{vw}{\widetilde{v}\widetilde{w}}.
\]
By the substitution%
\[
\widetilde{T}=\frac{T-T_{a}}{T_{1}},\qquad\widetilde{x}=\frac{x}{x_{1}}%
,\qquad\mbox{and}\qquad%
\widetilde{t}=\frac{t}{t_{1}}
\]
into the system (\ref{eq:heat2}-\ref{eq:fuel2}) with the coefficients%
\begin{equation}
A=T_{1}/t_{1},\quad B=T_{1},\quad C=\lambda/T_{1},\quad C_{S}=\beta%
/t_{1},\quad\mbox{and}\quad%
k=x_{1}^{2}/\left(  T_{1}^{3}t_{1}\right)%
\label{eq:inverse-scaling}%
\end{equation}
admits the scaled solution
\[
T(t,x)=T_{1}\widetilde{T}\left(  \frac{t}{t_{1}},\frac{x}{x_{1}}\right)
+T_{a},\quad\mbox{and}\quad
S(t,x)=\widetilde{S}\left(  \frac{t}{t_{1}},\frac{x}{x_{1}%
}\right)  ,
\]
which has the desired nondimensional properties as well as the correct maximal
temperature, width, and speed of a traveling combustion wave.

A dimensionless system similar to (\ref{eq:dimensionless-heat}%
-\ref{eq:dimensionless-fuel}) is studied in \cite{Weber-1997-CWG} in the case
$\lambda=0$, i.e., combustion insulated against heat loss.
\cite{Weber-1997-CWG} determine the speed of the traveling wave as a
function of $\beta$ numerically and by an asymptotic expansion and
observe that a traveling combustion wave exists only for small values of
$\beta$. By increasing $\beta$, the solution is periodic, then the period
doubles, and eventually the solution becomes chaotic. We have observed that
increasing $\beta$ has a similar effect for the equations
(\ref{eq:dimensionless-heat}-\ref{eq:dimensionless-fuel}) when $\lambda>0$.
Also, we have observed that a sustained combustion wave is possible only when
$\lambda$ is small enough. A systematic study of the properties of
(\ref{eq:dimensionless-heat}-\ref{eq:dimensionless-fuel}) for various values
of $\lambda$ and $\beta$ will be done elsewhere. For a dimensionless system
similar to ours, but without the temperature offset to force zero reaction at
ambient temperature, see \cite{Asensio-2002-WFM}.

\section{Data Assimilation}

\label{sec:enkf}

The goal in data assimilation on a fire model is a
filter that can effectively track the location of the fireline given data in
the form of temperature and remaining fuel at sample points inside of the
domain. The fire application is particularly troublesome for EnKFs. The
standard method for generating an initial ensemble is not sufficient for this
scenario. Namely, taking an initial guess at the model state (temperature and
fuel) and adding to it a smooth random field. Here, if the data indicates that
the fireline has shifted away from that of the ensemble, then the Kalman
Filter will generally ignore the data entirely due to the extraordinarily
small data likelihood. Clearly, such an initial ensemble does not properly
represent the prior uncertainty in the location of the ignition region, only
that of the temperature of ignition. In order to represent this uncertainty as
well, we have also perturbed the state variables by a spatial shift
\cite{Johns-2007-CEK}. However, this approach leads to further potential
problems for EnKF. Due to the relatively sharp temperature profile of the
fireline, the temperature at each grid point will tend to be close to that of
the stable ambient or burning temperatures. A similar situation occurs with
the fuel near the fireline as well. This is indicative of a strongly bimodal
or non-Gaussian prior distribution. Despite this violation of the underlying
assumptions of EnKF, we have found that it is possible to track large changes
in the fireline, as shown in the numerical results in Section \ref{sec:numerical-results}.

The Ensemble Kalman Filter (EnKF) is a Monte-Carlo implementation of the
Bayesian update problem:\ Given a probability distribution of the modeled
system (the prior, often called the `forecast' in geophysical sciences) and
data likelihood, the Bayes' theorem is used to obtain the probability
distribution with the data likelihood taken into account (the posterior or the
`analysis'). The Bayesian update is combined with advancing the model in time,
with the data incorporated from time to time. The original Kalman Filter
\cite{Kalman-1960-NAL} relies on the assumption that the probability
distributions are Gaussian (`the Gaussian assumption'), and provides algebraic formulas for the change of
the mean and covariance by the Bayesian update, and a formula for advancing
the covariance matrix in time provided the system is linear. However, this is
not possible computationally for high-dimensional systems. For this reason,
EnKFs were developed in \cite{Evensen-1994-SDA,Houtekamer-1998-DAE}. EnKFs
represent the distribution of the system state using an ensemble of simulations,
and replace the covariance matrix by the covariance matrix
of the ensemble. One advantage of EnKFs is that advancing the probability
distribution in time is achieved by simply advancing each member of the
ensemble. EnKFs, however, still rely on the Gaussian assumption, though they
are of course used in practice for nonlinear problems, where the Gaussian
assumption is not satisfied. Related filters attempting to relax the Gaussian
assumption in EnKF include
\cite{Anderson-1999-MCI,Beezley-2008-MEK,Bengtsson-2003-NFE,%
Mandel-2006-PEF,vanLeeuwen-2003-VMF}.

We use the EnKF following \cite{Burgers-1998-ASE,Evensen-2003-EKF}, with only
some minor differences. This filter involves randomization of data. For
filters without randomization of data, see
\cite{Anderson-2001-EAK,Evensen-2004-SSR,Tippett-2003-ESR}. The data
assimilation uses a collection of independent simulations, called an ensemble.
The ensemble filter consists of

\begin{enumerate}
\item generating an initial ensemble by random perturbations,

\item advancing each ensemble member in time until the time of the data, which
gives the so-called \emph{forecast ensemble}

\item modifying the ensemble members by injecting the data (the \emph{analysis
step}), which results in the so-called \emph{analysis ensemble}

\item continuing with step 2 to advance the ensemble in time again.
\end{enumerate}

We now consider the analysis step in more detail. We have the forecast
ensemble
\[
U^{f}=[u_{1}^{f},\ldots,u_{N}^{f}]=[u_{i}^{f}]
\]
where each $u_{i}^f$ is a column vector of dimension $n$, which contains the whole
simulation state (in our case, the vector of the values of $T$ and $S$ at mesh
nodes). Thus, $U^{f}$ is a matrix of dimension $n$ by $N$. The superscript
$^{f}$ stands for \textquotedblleft forecast\textquotedblright. The data is
given as a measurement vector $d$ of dimension $m$ and data error covariance
matrix $R$ of dimension $m$ by $m$. The correspondence of the data and the
simulation states is given by an \emph{observation function} $h(u)$ that
creates synthetic data, that is, what the data would have been if the
simulation and the measurements were exact. We assume that $h$ is linear,
\begin{equation}
h\left(  u\right)  =Hu. \label{eq:h}%
\end{equation}
Using the notation that $\propto$ means proportional and
$N(x,M)$ represents a~normal distribution with mean $x$ and covariance
matrix $M$,
the observation being assimilated is%
\begin{equation}
Hu-d\sim N\left(  0,R\right)  \label{eq:data-error}%
\end{equation}
for some matrix $H$. The observation function defines the data likelihood (the
probability density of $d$ given the state $u$). Assuming the data error is
normally distributed, the data likelihood is%
\[
p\left(  d|u\right)  \propto e^{-\frac{1}{2}\left(  h\left(  u\right)
-d\right)  R^{-1}\left(  h\left(  u\right)  -d\right)  }.
\]
The forecast ensemble $U^{f}$ is
considered a sample from the prior distribution $p(u)$, and the EnKF
strives to create an analysis ensemble that is a sample from the
posterior distribution $p\left(  u|d\right)  $, which is the probability
distribution of $u$ after the data has been injected. From Bayes' theorem of
probability theory, we have%
\begin{equation}
p\left(  u|d\right)  \propto p\left(  d|u\right)  p(u) \label{eq:bayes},%
\end{equation}
\newline \cite[eq. (20)]{Evensen-2003-EKF}. If $p\sim N(u^{f},Q^{f})$,
then it is known that the posterior is also normally
distributed with mean%
\begin{equation}
u^{a}=u^{f}+K\left(  d-Hu^{f}\right)  , \label{eq:Kalman-update}%
\end{equation}
where $K$ is the Kalman gain matrix,%
\begin{equation}
K=Q^{f}H^{\mathrm{T}}\left(  HQ^{f}H^{\mathrm{T}}+R\right)  ^{-1}.
\label{eq:Kalman-gain}%
\end{equation}

The EnKF is based on applying a version of the Kalman update
(\ref{eq:Kalman-update})\ to each forecast ensemble member $u_{i}^{f}$ to
yield the analysis ensemble member $u_{i}^{a}$. For this update, the data
vector $d$ in (\ref{eq:Kalman-update}) is replaced by a randomly perturbed
vector%
\begin{equation}
d_{j}=d+v_{j},\quad v_{j}\sim N\left(  0,R\right)  . \label{eq:data-pert}%
\end{equation}
Let $\bar{u}^{f}$ be the mean of the forecast ensemble,
\begin{equation}
\bar{u}^{f}=\frac{1}{N}\sum_{i=1}^{N}u_{i}^{f}. \label{eq:mean}%
\end{equation}
The unknown covariance matrix $Q^{f}$ in (\ref{eq:Kalman-gain}) is
replaced by the covariance matrix $C^{f}$ of the forecast ensemble $U^{f}$,
\begin{equation}
C^{{}}=\frac{1}{N-1}AA^{T},\quad A=[u_{1}^{f}-\bar{u}^{f},\ldots,u_{N}%
^{f}-\bar{u}^{f}]. \label{eq:covariance}%
\end{equation}
Define%
\[
D=\left[  d+v_{1},\ldots,d+v_{N}\right]
\]
as the matrix formed from the randomly perturbed data vectors and%
\[
U^{a}=[u_{1}^{a},\ldots,u_{N}^{a}]=[u_{i}^{a}]
\]
as the analysis ensemble.
This gives the EnKF formula,
\begin{equation}
U^{a}=U^{f}+CH^{T}\left(  HCH^{T}+R\right)  ^{-1}(D-HU^{f}). \label{eq:enkf}%
\end{equation}

See \cite[eq. (20)]{Evensen-2003-EKF} for details. The only difference between
(\ref{eq:enkf}) and \cite[eq. (20)]{Evensen-2003-EKF} is that we use the
covariance matrix $R$ of the measurement error (which is assumed to be known anyway)
rather than the sample
covariance matrix of the randomized data. Since $R$ is in practice
positive definite, there is no difficulty with the inverse in (\ref{eq:enkf}).
The matrix $R$ is known from (\ref{eq:data-error}). In \cite[eq.
(20)]{Evensen-2003-EKF}, the sample covariance matrix, called $R_{e}$, of the
perturbed data is used in place of $R$. For large $n$, the matrix $R_{e}$ is
singular and then a more expensive pseudoinverse using
an eigenvalue decomposition
\cite[eq. (56)]{Evensen-2003-EKF} must be used.
Alternately, the data perturbations
have to be chosen in a special way \cite[eq. (57)]{Evensen-2003-EKF}.

\subsection{EnKF implementation}

We have used the ensemble update (\ref{eq:enkf}) with the inverse computed by
an application of the Sherman-\allowbreak Morrison-\allowbreak Wood\-bury
formula~\cite{Hager-1989-UIM}%
\begin{align}
\lefteqn{\left(  HCH^{T}+R\right)  ^{-1}  =\left(  R+\frac{1}{N-1}HA\left(
HA\right)  ^{T}\right)  ^{-1}\nonumber =}\\
&R^{-1}\left[  I-\frac{1}{N-1}\left(  HA\right)  \left(  I+\left(
HA\right)  ^{T}R^{-1}\frac{1}{N-1}\left(  HA\right)  \right)  ^{-1}\left(
HA\right)  ^{T}R^{-1}\right]  . \label{eq:P-inv}%
\end{align}
%% \begin{align}
%% \left(  HCH^{T}+R\right)  ^{-1}  &  =\left(  R+\frac{1}{N-1}HA\left(
%% HA\right)  ^{T}\right)  ^{-1}\nonumber\\
%% &  =R^{-1}\left[  I-\frac{1}{N-1}\left(  HA\right)  \left(  I+\left(
%% HA\right)  ^{T}R^{-1}\frac{1}{N-1}\left(  HA\right)  \right)  ^{-1}\left(
%% HA\right)  ^{T}R^{-1}\right]  . \label{eq:P-inv}%
%% \end{align}
This formula is advantageous when the data error covariance matrix $R$ is of a
special form such that left and right multiplications by $R^{-1}$ can be
computed inexpensively. In particular, when the data errors are uncorrelated,
which is usually the case in practice and the case here, the matrix $R$ is
diagonal. Then the EnkF formula (\ref{eq:enkf}) with (\ref{eq:P-inv}) costs
$O\left(  N^{3}+mN^{2}+nN^{2}\right)  $ operations, which is suitable both for
a large number $n$ of the degrees of freedom and a large number $m$ of data
points. Also, (\ref{eq:enkf})
can be implemented without forming the observation matrix $H$ explicitly by
only evaluating the observation function $h$ using%
\begin{align*}
\left[  HA\right]  _{i}  &  =Hu_{i}^{f}-H\frac{1}{N}\sum_{j=1}^{N}u_{i}%
^{f}=h\left(  u_{i}^{f}\right)  -\frac{1}{N}\sum_{j=1}^{N}h\left(  u_{i}%
^{f}\right)  ,\\
d-Hu_{i}^{f}  &  =d-h\left(  u_{i}^{f}\right)  .
\end{align*}
See \cite{Mandel-2006-EIE} for further details.

The ensemble filter formulas are operations on full matrices, and they were
implemented in a distributed parallel environment using MPI and ScaLAPACK.
EnKF is naturally parallel:\ each ensemble member can be advanced in time
independently. The linear algebra in the Bayesian update step links the
ensemble members together.

\subsection{Regularization}

EnKF produces the analysis ensemble in the span of the forecast ensemble. This
results in nonphysical states especially if the states in the span are far
away from the data. For cheap numerical methods and a highly nonlinear
problem, EnKF can easily knock the state out of the stability region. In order
to ease this problem, we add an independent observation%
\[
\bigtriangledown u^{a}-\bigtriangledown\bar{u}^{f}\sim N(0,D),
\]
where $\bigtriangledown$ is the spatial gradient, computed by finite
differences. This is easily implemented by running the EnKF\ formulas a second
time. In practice, this matrix, $D$, is of the form $\rho I$, where $\rho$ is
a regularization parameter. This technique prevents large, nonphysical
gradients in the analysis ensemble. See \cite{Johns-2007-CEK} for further details.

\section{Numerical results}

\label{sec:numerical-results}

\subsection{Calibration of coefficients in one dimension}

\begin{figure}[ptb]
\begin{center}
\includegraphics[width=3in]{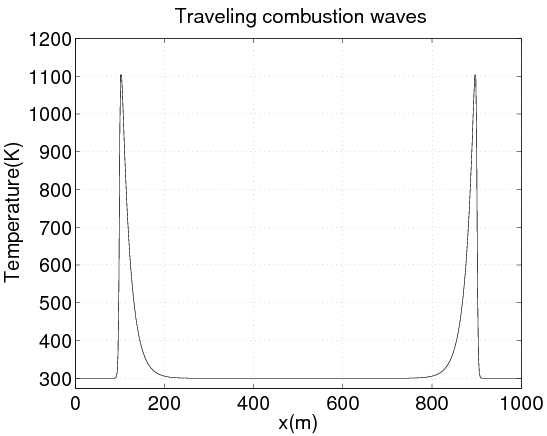}
\end{center}
\caption{Temperature profile of traveling wave. The wave moved about 398{m}
from the initial position in 2300{s}.}%
\label{fig:tw}%
\end{figure}

\begin{figure}[ptb]
\begin{center}
\includegraphics[width=3in]{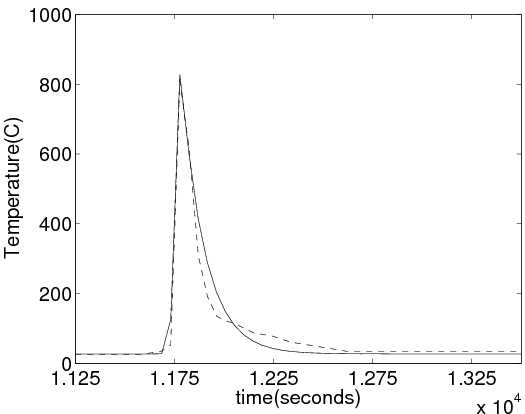}
\end{center}
\par
\bigskip\caption{Time-temperature profile (dotted line) measured in a grass
wildland fire at a fixed sensor location, digitized from
\cite{Kremens-2003-MTT}, and a computed profile (solid line) from simulation.}%
\label{fig:tp}%
\end{figure}

We have found initial values $B=5.5849\times10^{4} K$, and $C=5.9739\times
10^{-4} K^{-1}$ from the values $T_{i}=670 K$ and $T_{c}=1200 K$ using
(\ref{eq:BC}), and then the value $A=1.5217\times10^{1} K s^{-1}$ from
(\ref{eq:A}), using the value $t_{c}=110 s$\/\ from \cite{Kremens-2003-MTT}. Not
every initial condition gives rise to a traveling combustion wave
\cite{Mercer-1996-CPW}. Inspired by \cite{Mercer-1997-CWT}, we use
an initial condition of the form
$T(x,t_{0})=T_{c}e^{-\left( x-x_{0}\right) ^{2}/\sigma^{2}}+T_{a}$,
where $x_{0}$ is in the center of
the interval and $\sigma=10\sqrt{2}m$. This initial condition is smooth. Thus,
it does not excite possible numerical artifacts. It has numerically local
support and for a modest $\sigma$ provides ignition sufficient to
develop into two sustained combustion waves traveling from the center. We have
then found empirically suitable values of $k$ and $C_{S}$ that result in
traveling combustion waves, computed the nondimensional coefficients $\lambda$
and $\beta$, adjusted them using Fig.~\ref{fig:tp} (dotted line), and scaled
using (\ref{eq:inverse-scaling}) to match the maximal temperature and the
width of the wave in Fig.~\ref{fig:tp} (dotted line), and the speed of the
traveling combustion wave, $0.17m/s$, from \cite{Kremens-2003-MTT}. There was
a small amount of wind in \cite{Kremens-2003-MTT}, however we have not
considered the wind here. The resulting coefficients are $k=2.1360\times
10^{-1}m^{2}s^{-1}K^{-3}$, $A=1.8793\times10^{2}Ks^{-1}$, $B=5.5849\times
10^{2}K$, $C=4.8372\times10^{-5}K^{-1}$, and $C_{S}=1.6250\times10^{-1}s^{-1}%
$. The corresponding nondimensional coefficients were $\lambda=2.7000\times
10^{-2}$ and $\beta=0.4829$. The computed traveling combustion wave
(see Figs.~\ref{fig:ra}, \ref{fig:tw} and \ref{fig:tp}, solid line) is a
reasonable match with the observation (see Fig.~\ref{fig:tp}, dotted line). The
trailing edge of the computed temperature profile (see Fig.~\ref{fig:tp}, solid
line) was not so well matched but this model is quite limited and other
matches reported in the literature are similar \cite{Balbi-1999-DMF}. The real
data looks like the superposition of two exponential decay modes, possibly the
fast one from cooling and the long one from the heat stored in water in the
ground.

We have also noted that when the ratio $A/C_{S}$ increases the temperature in
the traveling combustion wave increases, increasing the thermal diffusivity
coefficient $k$ increases the width and speed of the combustion wave and that
the maximum temperature in the traveling wave decreases if $C_{S}$ increases.
Sufficiently small value of $C_{S}$ is needed for sustained combustion. We
have noted that the numerical solution by finite differences becomes unstable
when the ratio $k/h$, where $h$ is the mesh size, is too small.

\subsection{Numerical results in two dimensions}

We have implemented the fire model in two dimensions by central finite
differences in space. The mesh size was $250$ by $250$ and the mesh step was
$2 m$. We have used the explicit Euler method with time step $1 s$. The
initial conditions were given by the ambient temperature $T_{a}=300 K$
everywhere except in a $50m{\times}50m$
square ignition region which was ignited
by elevating the temperature to $1200K$. The mass fraction of the fuel was
initialized to be one everywhere except for a $25m$ fuel break in the center
of the domain. Then, at each grid point, the fuel was shifted by a random
number in $[-0.3,0.3]$.
This is intended to simulate a natural uniform fuel
supply and a road as a fuel break. The Neumann boundary conditions were
specified on all boundaries with no ambient wind across the domain.

The initial ensemble was generated by perturbing the temperature profile of
what we call the \textit{comparison solution} $T_{0}$ utilizing smooth
random fields in the following form:%

\begin{equation}
\tilde{u}=\sum_{n=1}^{d}\frac{v_{n}}{1+n^{2\alpha}}e_{n},\quad v_{n}\sim
N\left(  0,1\right)  , \label{eq:smooth_field}%
\end{equation}
where $\alpha$ is the order of smoothness of the random field, and $\left\{
e_{n}\right\}  _{n=1}^{d}$ is the Fourier sine basis, ensuring that $\tilde
{u}$ is real valued \cite{Evensen-1994-SDA,Ruan-1998-EMR}. This process can be
understood as a finite dimensional version of sampling out of normal
distribution on an infinite dimensional space of smooth functions, the Sobolev
space $H_{0}^{\alpha}\left(  \Omega\right)  $ \cite{Mandel-2006-PEF}. For
integer $\alpha$, this is the space of functions with square integrable
derivatives of order $\alpha$ and zero traces on the boundary.

\begin{figure}[ptb]
\begin{center}%
\begin{tabular}
[c]{cc}%
\includegraphics[width=2.5in]{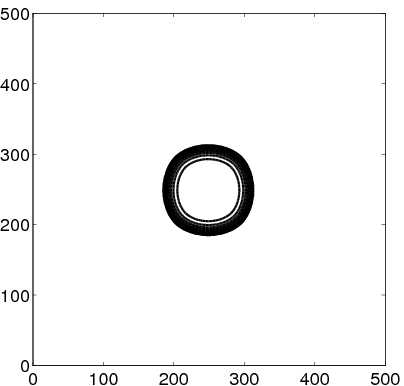} &
\includegraphics[width=2.5in]{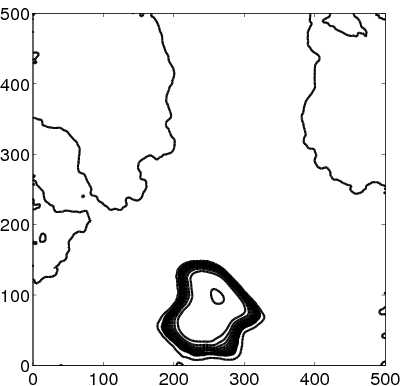}\\
(a) & (b)
\end{tabular}
\end{center}
\caption{Contour plots with $100K$ between contour lines. %
(a) The temperature profile of a circular ignition region in the center of %
domain. (b) The same profile randomly perturbed in magnitude by (\ref{eq:T-pert}) with
 $c_{T}=5$ and spatially by (\ref{eq:xy-pert}) with $c_{x}=c_{y}=100$.}%
\label{fig:da_pert}%
\end{figure}

A preliminary ensemble was generated by adding a smooth random field to each
state variable of the comparison solution. For example, the temperature of
$k^{th}$ ensemble member is given by
\begin{equation}
\hat{T}_{k}=T_{0}+c_{T}\tilde{u}_{k},
\label{eq:T-pert}
\end{equation}
where the scalar $c_{T}$ controls the magnitude of this perturbation. Finally,
the preliminary ensemble was moved spatially in both $x$ and $y$ directions
by
\begin{equation}
T_{k}(x,y)=\hat{T}_{k}(x+c_{x}\tilde{u}_{i_{k}}(x,y),y+c_{y}\tilde{u}_{j_{k}%
}(x,y)).
\label{eq:xy-pert}
\end{equation}
Here, $c_{x}$ and $c_{y}$ control the magnitude of the shift in each
coordinate, bilinear interpolation is used to determine $\hat{T}$ on off-grid
points, and the temperature outside of the computational domain is assumed to
be at the ambient temperature. The given simulation was run with the
initialization parameters $c_{T}=5$ and $c_{x}=c_{y}=150$.
Fig.~\ref{fig:da_pert} shows the effect of these perturbations on a simple
circular fire line on the center of the domain.

In each analysis cycle, the solution was advanced by $100s$, and then the data
was injected. The data was created artificially by sampling the temperature
and fuel of one fixed solution, called the \textit{reference solution}, every
$10m$ (or 5 grid points) across the domain. The data covariance matrix was
taken to be diagonal with a variance $10$ for each sample, and the
regularization was used with regularization parameter $\rho=750$. The
reference solution was created in the same manner as the ensemble with an
ignition region located $100m$ away in one direction. This discrepancy is
intended to demonstrate the power of EnKF to attract the ensemble to the
truth. After each analysis cycle, the ensemble was further perturbed by $5\%$
magnitude of the initial perturbation to assure sufficient ensemble spread for
future assimilations.

\begin{figure}[ptb]
\begin{center}%
\begin{tabular}
[c]{cc}%
\includegraphics[width=2.5in]{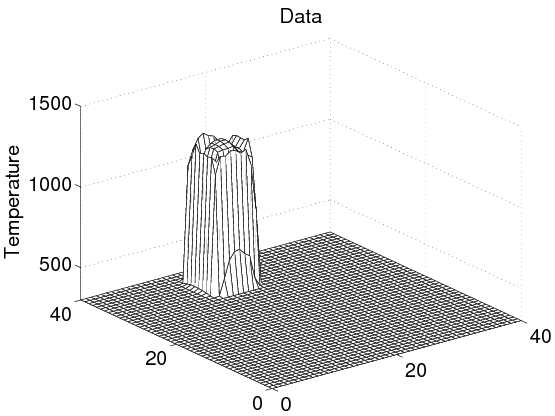} &
\includegraphics[width=2.5in]{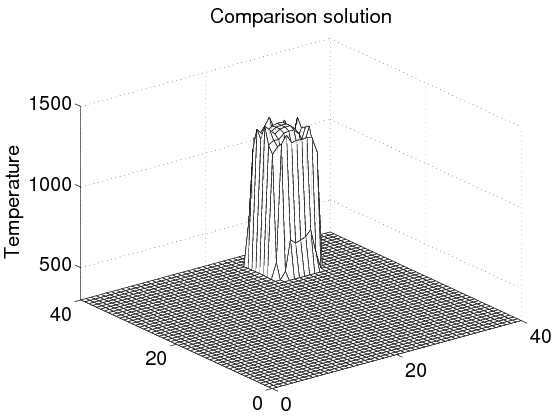}\\
(a) & (b)
\end{tabular}
\end{center}
\caption{Temperature profiles representing (a) the data (reference solution,
taken as the truth) and (b) an unperturbed ensemble member comparison solution
$100s$ after initialization.}%
\label{fig:da_refcomp}%
\end{figure}

\begin{figure}[ptb]
\begin{center}%
\begin{tabular}
[c]{cc}%
\includegraphics[width=2.5in]{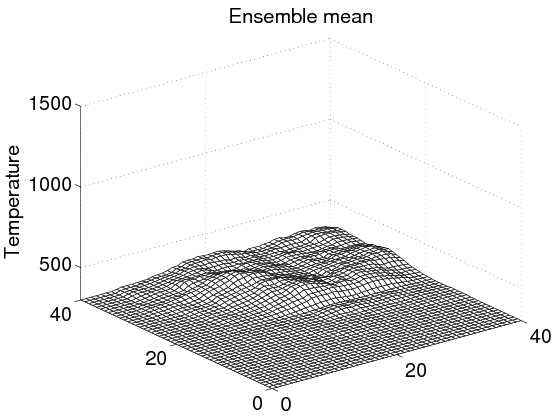} &
\includegraphics[width=2.5in]{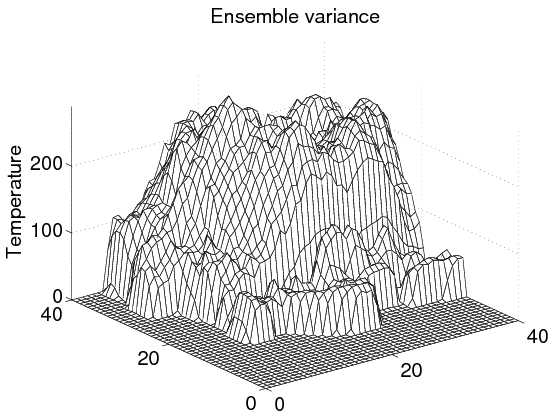}\\
(a) & (b)
\end{tabular}
\end{center}
\caption{After advancing the solution in time by $100s$ before any %
data assimilations: %
Pointwise prior ensemble (a) mean and (b) variance.}%
\label{fig:da_ens1}%
\end{figure}

\begin{figure}[ptb]
\begin{center}%
\begin{tabular}
[c]{cc}%
\includegraphics[width=2.5in]{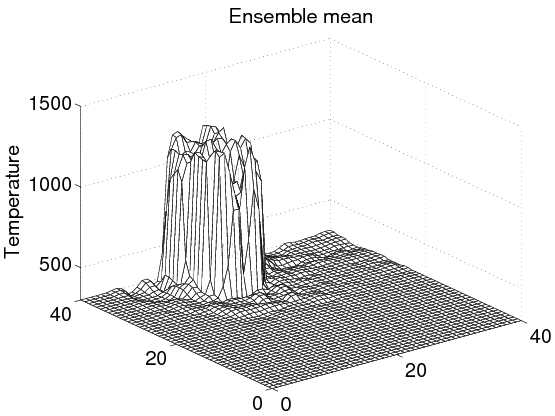} &
\includegraphics[width=2.5in]{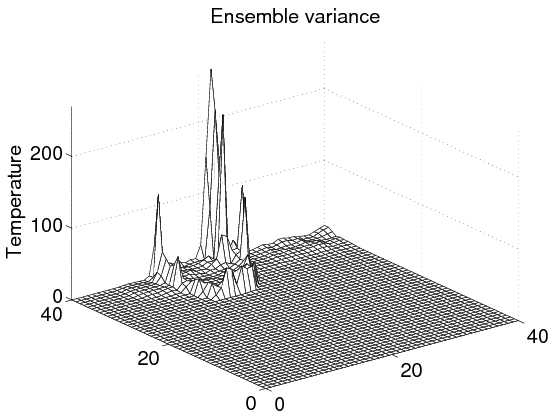}\\
(a) & (b)
\end{tabular}
\end{center}
\caption{After advancing the solution in time by $100s$ and performing a %
single data assimilation: %
Pointwise posterior ensemble (a) mean and (b) variance.}%
\label{fig:da_ens2}%
\end{figure}

\begin{figure}[ptb]
\begin{center}%
\begin{tabular}
[c]{cc}%
\includegraphics[width=2.5in]{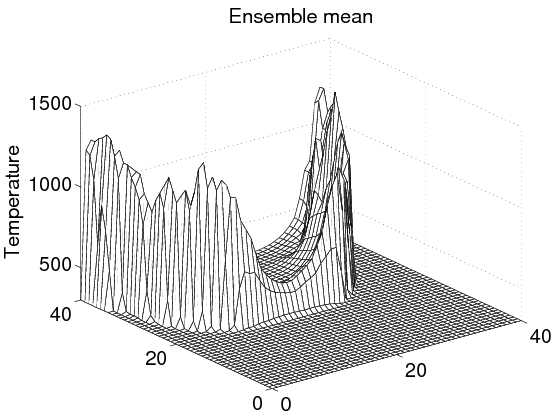} &
\includegraphics[width=2.5in]{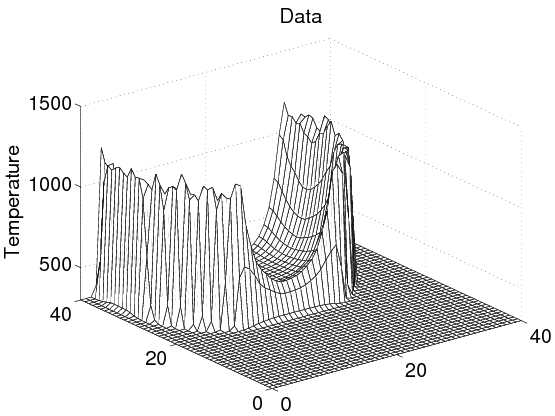}\\
(a) & (b)
\end{tabular}
\end{center}
\caption{After $10$ analysis cycles with a $100s$ time update per cycle, (a)
the ensemble mean compared to (b) the reference solution (the data).}%
\label{fig:da_ens3}%
\end{figure}

\begin{figure}[ptb]
\begin{center}
\includegraphics[width=2.5in]{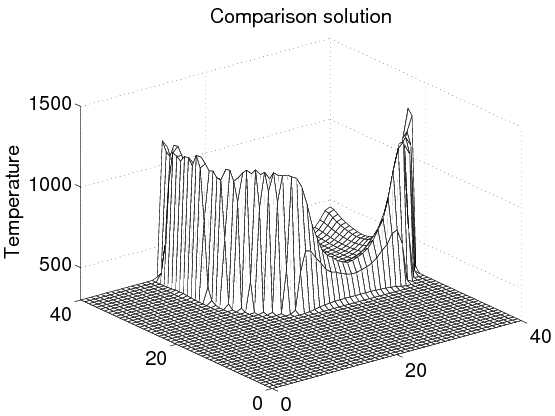}
\end{center}
\caption{Comparison solution advanced to $1000s$. This is what the solution in
Fig.~\ref{fig:da_ens3} (a) would be without data assimilation.}%
\label{fig:da_comp3}%
\end{figure}

Fig.~\ref{fig:da_refcomp} shows the reference and comparison solution $100s$
after initialization, at the end of the first cycle.
Fig.~\ref{fig:da_ens1} shows the ensemble mean and variance at the same time
prior to performing an assimilation. Fig.~\ref{fig:da_ens2} shows the ensemble
after applying the first assimilation. The analysis cycle was repeated $10$
times with the results shown in Figs.~\ref{fig:da_ens3}
and~\ref{fig:da_comp3}. These figures show a remarkable agreement of the
ensemble mean with the reference solution, even if the simulation ensemble was
ignited intentionally far away from the reference ignition region. However, it
should be noted that different runs of this \textit{stochastic} algorithm
produce different results. Sometimes the ensemble is attracted to the
reference solution, and sometimes not, depending on if there exists a good
match to the data in the span of a fairly small ensemble.

\section{Conclusion}

\label{sec:conclusions}

A simple model based on two coupled PDEs can reproduce the time-temperature
curve recorded as a wildfire burns over a sensor, which is a measurable
feature of fire behavior. By separating the parameters that determine the
qualitative properties of the solution from the parameters that determine the
temperature, time, and space scales, we were able to identify the parameters
of the model from actual wildfire observations. Assimilation of data into a
wildfire simulation poses a particular challenge because the combustion
region is quite thin. We have shown that a version of the Ensemble Kalman
filter is able to assimilate data into a wildfire simulation successfully. The
filter uses penalization of nonphysical solution and perturbations by smooth
random transformations of the spatial domain, in addition to standard smooth
additive perturbation.

\section{Acknowledgments}

This material is based upon work supported by the National Science Foundation
(NSF) under grants CNS-0325314, CNS-0324989, CNS-0324876, CNS-0324910,
CNS-0720454, CNS-0719641, CNS-0719626,
EIA-0219627, ACI-0305466, OISE-0405349, CNS-0540178, DMS-0610039, and by a
National Center for Atmospheric Research (NCAR) Faculty Fellowship. Computer
time on IBM BG/L was provided in part by NSF MRI Grants CNS-0421498,
CNS-0420873, CNS-0420985, NSF sponsorship of the NCAR, the University of
Colorado, and a grant from the IBM Shared University Research (SUR) program.
TeraGrid computer time was provided by the NSF through TeraGrid resources at
the National Center for Supercomputing Applications, the San Diego
Supercomputing Center, the Texas Advanced Computing Center, and the Pittsburgh
Supercomputing Center. The authors would like to thank Bob Kremens for
suggesting an explanation for the behavior of the measured temperature in
Fig.~\ref{fig:tp}. The authors would like to thank an anonymous referee for
useful comments that contributed to improving this paper.

%\nocite{Hacker-2005-EKF}

\bibliographystyle{elsart-num-sort}
\bibliography{../../bibliography/dddas-jm}

\def\cprime{$'$} \def\cprime{$'$}
\begin{thebibliography}{100}
\expandafter\ifx\csname url\endcsname\relax
  \def\url#1{\texttt{#1}}\fi
\expandafter\ifx\csname urlprefix\endcsname\relax\def\urlprefix{URL }\fi

\bibitem{Albini-1986-WFS}
F.~A. Albini, A model for fire spread in wildland fuels by radiation - {A}
  model including fuel cooling by convection, Combustion and Science Technology
  45 (1986) 101--113.

\bibitem{Albini-1994-PBS}
F.~A. Albini, {PROGRAM BURNUP}: {A} simulation model of the burning of large
  woody natural fuels, final Report on Research Grant INT-92754-GR by U.S.F.S.
  to Montana State Univ., Mechanical Engineering Dept. (1994).

\bibitem{Albini-1997-ICL}
F.~A. Albini, E.~D. Reinhardt, Improved calibration of a large fuel burnout
  model, Int. J. Wildland Fire 7 (1997) 21--28.

\bibitem{Anderson-2001-EAK}
J.~L. Anderson, An ensemble adjustment {K}alman filter for data assimilation,
  Monthly Weather Review 129 (1999) 2884--2903.

\bibitem{Anderson-1999-MCI}
J.~L. Anderson, S.~L. Anderson, A {M}onte {Ca}rlo implementation of the
  nonlinear filtering problem to produce ensemble assimilations and forecasts,
  Monthly Weather Review 127 (1999) 2741--2758.

\bibitem{Asensio-2000-TEE}
M.~I. Asensio, L.~Ferragut, Total error estimates of mixed finite element
  methods for nonlinear reaction-diffusion equations, Neural Parallel Sci.
  Comput. 8~(2) (2000) 169--189.

\bibitem{Asensio-2002-WFM}
M.~I. Asensio, L.~Ferragut, On a wildland fire model with radiation, Int. J.
  Numer. Meth. Engrg. 54 (2002) 137--157.

\bibitem{Asensio-2004-RBM}
M.~I. Asensio, A.~Russo, G.~Sangalli, The residual-free bubble numerical method
  with quadratic elements, Mathematical Models and Methods in Applied Sciences
  14~(5) (2004) 641--661.

\bibitem{Baines-1990-PMP}
P.~G. Baines, Physical mechanisms for the propagation of surface fires,
  Mathematical and Computer Modelling 13 (1990) 83--94.

\bibitem{Balbi-1999-DMF}
J.~H. Balbi, P.~A. Santoni, J.~L. Dupuy, Dynamic modelling of fire spread
  across a fuel bed, International J. of Wildland Fire 9~(4) (1999) 275--284.

\bibitem{Beezley-2008-MEK}
J.~D. Beezley, J.~Mandel, Morphing ensemble {K}alman filters, Tellus 60A (2008)
  131--140.

\bibitem{Bengtsson-2003-NFE}
T.~Bengtsson, C.~Snyder, D.~Nychka, Toward a nonlinear ensemble filter for high
  dimensional systems, Journal of Geophysical Research - Atmospheres 108(D24)
  (2003) STS 2--1--10.

\bibitem{Berestycki-1991-MIC}
H.~Berestycki, B.~Larrouturou, J.-M. Roquejoffre, Mathematical investigation of
  the cold boundary difficulty in flame propagation theory, in: P.~C. Fife,
  A.~Li{\~ n}\'{a}n, F.~Williams (eds.), Dynamical issues in combustion theory
  (Minneapolis, MN, 1989), vol.~35 of IMA Vol. Math. Appl., Springer, New York,
  1991, pp. 37--61.

\bibitem{Burgers-1998-ASE}
G.~Burgers, P.~J. van Leeuwen, G.~Evensen, Analysis scheme in the ensemble
  {K}alman filter, Monthly Weather Review 126 (1998) 1719--1724.

\bibitem{Campos-2004-RDP}
D.~Campos, J.~E. Llebot, J.~Fort, Reaction-diffusion pulses: a combustion
  model, J. Phys. A: Math. Gen. 37 (2004) 6609--6621.

\bibitem{Carey-1995-LFE}
G.~F. Carey, Y.~Shen, Least-squares finite element approximation of {F}isher's
  reaction-diffusion equation, Numer. Methods Partial Differential Equations
  11~(2) (1995) 175--186.

\bibitem{Chen-1992-GPIs}
X.~Chen, Generation and propagation of interfaces in reaction-diffusion
  systems, Trans. Amer. Math. Soc. 334~(2) (1992) 877--913.

\bibitem{Chen-2007-AVP}
Y.~Chen, C.~Snyder, Assimilating vortex position with an ensemble {K}alman
  filter, Monthly Weather Review 135 (2007) 1828--1845.

\bibitem{Clark-2004-DCA}
T.~L. Clark, J.~Coen, D.~Latham, Description of a coupled atmosphere-fire
  model, Intl. J. Wildland Fire 13 (2004) 49--64.

\bibitem{Clark-1996-CAF}
T.~L. Clark, M.~A. Jenkins, J.~Coen, D.~Packham, A coupled atmospheric-fire
  model: {C}onvective feedback on fire line dynamics, J. Appl. Meteor 35 (1996)
  875--901.

\bibitem{Class-2003-UMF}
A.~G. Class, B.~J. Matkowsky, A.~Y. Klimenko, A unified model of flames as
  gasdynamic discontinuities, J. Fluid Mech. 491 (2003) 11--49.

\bibitem{Codina-1998-CFE}
R.~Codina, Comparison of some finite element methods for solving the
  diffusion-convection-reaction equation, Comput. Methods Appl. Mech. Engrg.
  156~(1-4) (1998) 185--210.

\bibitem{Darema-2004-DDD}
F.~Darema, Dynamic data driven applications systems: {A} new paradigm for
  application simulations and measurements, in: M.~Bubak, G.~D. van Albada,
  P.~M.~A. Sloot, J.~J. Dongarra (eds.), Computational Science-ICCS 2004: 4th
  International Conference, vol. 3038 of Lecture Notes in Computer Science,
  Springer, 2004, pp. 662--669.

\bibitem{Dold-2003-HOE}
J.~W. Dold, R.~W. Thatcher, A.~A. Shah, High order effects in one step reaction
  sheet jump conditions for premixed flames, Combust. Theory Model. 7~(1)
  (2003) 109--127.

\bibitem{Douglas-2006-DVW}
C.~C. Douglas, J.~D. Beezley, J.~Coen, D.~Li, W.~Li, A.~K. Mandel, J.~Mandel,
  G.~Qin, A.~Vodacek, Demonstrating the validity of a wildfire {DDDAS}, in:
  V.~N. Alexandrov, D.~G. van Albada, P.~M.~A. Sloot, J.~Dongarra (eds.),
  Computational Science – ICCS 2006: 6th International Conference, Reading, UK,
  May 28-31, 2006, Proceedings, Part III, vol. 3993 of Lecture Notes in
  Computer Science, Springer, Berlin/Heidelberg, 2006, pp. 522--529.

\bibitem{Dupuy-2000-TTR}
J.-L. Dupuy, Testing two radiative physical models for fire spread through
  porous forest fuel beds, Combustion Science and Technology 155~(1) (2000)
  149--180.

\bibitem{Dupuy-1999-FSP}
J.-L. Dupuy, M.~Larini, Fire spread through a porous forest fuel bed: A
  radiative and convective model including fire-induced flow effects,
  International J. of Wildland Fire 9~(3) (1999) 155--172.

\bibitem{Dupuy-2005-NSC}
J.-L. Dupuy, D.~Morvan, Numerical study of a crown fire spreading toward a fuel
  break using a multiphase physical model, International Journal of Wildland
  Fire 14 (2005) 141–--151.

\bibitem{Emara-2002-NMS}
H.~E. Emara-Shabaik, Y.~A. Khulief, I.~Hussaini, A non-linear multiple-model
  state estimation scheme for pipeline leak detection and isolation,
  Proceedings of the Institution of Mechanical Engineers, Part I: Journal of
  Systems and Control Engineering 216 (2002) 497--512.

\bibitem{Ern-1995-DCM}
A.~Ern, C.~C. Douglas, M.~D. Smooke, Detailed chemistry modeling of laminar
  diffusion flames on parallel computers, The International Journal of
  Supercomputer Applications and High Performance Computing 9 (1995) 167--186.

\bibitem{Evensen-1994-SDA}
G.~Evensen, Sequential data assimilation with nonlinear quasi-geostrophic model
  using {M}onte {C}arlo methods to forecast error statistics, Journal of
  Geophysical Research 99 (C5)~(10) (1994) 143--162.

\bibitem{Evensen-2003-EKF}
G.~Evensen, The ensemble {K}alman filter: {T}heoretical formulation and
  practical implementation, Ocean Dynamics 53 (2003) 343--367.

\bibitem{Evensen-2004-SSR}
G.~Evensen, Sampling strategies and square root analysis schemes for the
  {EnKF}, Ocean Dynamics 54 (2004) 539--560.

\bibitem{Ferragut-2002-MFE}
L.~Ferragut, I.~Asensio, Mixed finite element methods for a class of nonlinear
  reaction diffusion problems, Neural Parallel Sci. Comput. 10~(1) (2002)
  91--112.

\bibitem{Fife-1988-DIL}
P.~C. Fife, Dynamics of internal layers and diffusive interfaces, vol.~53 of
  CBMS-NSF Regional Conference Series in Applied Mathematics, Society for
  Industrial and Applied Mathematics (SIAM), Philadelphia, PA, 1988.

\bibitem{Franca-2005-TMS}
L.~P. Franca, A.~L. Madureira, F.~Valentin, Towards multiscale functions:
  enriching finite element spaces with local but not bubble-like functions,
  Comput. Methods Appl. Mech. Engrg. 194~(27-29) (2005) 3006--3021.

\bibitem{Franca-2006-EFM}
L.~P. Franca, J.~V.~A. Ramalho, F.~Valentin, Enriched finite element methods
  for unsteady reaction-diffusion problems, Communications in Numerical Methods
  in Engineering 22 (2006) 519--526.

\bibitem{Frandsen-1971-SFP}
W.~H. Frandsen, Fire spread through porous fuels from conservation of energy,
  Combustion and Flame 16 (1971) 9--16.

\bibitem{Frank-Kamenetskii-1955-DHE}
D.~A. Frank-Kamenetskii, Diffusion and heat exchange in chemical kinetics,
  Princeton University Press, 1955.

\bibitem{Gazdag-1974-NSF}
J.~Gazdag, J.~Canosa, Numerical solution of {F}isher's equation, J. Appl.
  Probability 11 (1974) 445--457.

\bibitem{Gilding-2004-TWN}
B.~H. Gilding, R.~Kersner, Travelling waves in nonlinear diffusion-convection
  reaction, Progress in Nonlinear Differential Equations and their
  Applications, 60, Birkh\"auser Verlag, Basel, 2004.

\bibitem{Giroud-1998-DAT}
F.~Giroud, J.~Margerit, C.~Picard, O.~S\'{e}ro-Guillaume, Data assimilation:
  {T}he need for a protocole, in: D.~X. Viegas (ed.), Forest Fire Research:
  Proceedings 3rd International Conference on Forest Fire Research and 14th
  Conference on Fire and Forest Meteorology, Louso, Coimbra, Portugal, 16--18
  November, 1998, vol.~1, Associa\c{c}\~{a}o para o Desenvolvimento da
  Aerodinamica Industrial, 1998, pp. 583--598.

\bibitem{Grishin-1996-GMM}
A.~M. Grishin, General mathematical model for forest fires and its
  applications, Combustion Explosion and Shock Waves 32 (1996) 503--519.

\bibitem{Grishin-2002-MMS}
A.~M. Grishin, O.~V. Shipulina, Mathematical model for spread of crown fires in
  homogeneous forests and along openings, Combustion Explosion and Shock Waves
  38 (2002) 622--632.

\bibitem{Gubernov-2003-EFS}
V.~Gubernov, G.~N. Mercer, H.~S. Sidhu, R.~O. Weber, Evans function stability
  of combustion waves, SIAM J. Appl. Math. 63~(4) (2003) 1259--1275.

\bibitem{Gubernov-2004-EFS}
V.~V. Gubernov, G.~N. Mercer, H.~S. Sidhu, R.~O. Weber, Evans function
  stability of non-adiabatic combustion waves, Proc. R. Soc. Lond. Ser. A Math.
  Phys. Eng. Sci. 460~(2048) (2004) 2415--2435.

\bibitem{Hager-1989-UIM}
W.~W. Hager, Updating the inverse of a matrix, SIAM Rev. 31~(2) (1989)
  221--239.

\bibitem{Houtekamer-1998-DAE}
P.~Houtekamer, H.~L. Mitchell, Data assimilation using an ensemble {K}alman
  filter technique, Monthly Weather Review 126~(3) (1998) 796--811.

\bibitem{Infeld-2000-NWS}
E.~Infeld, G.~Rowlands, Nonlinear waves, solitons and chaos, 2nd ed., Cambridge
  University Press, Cambridge, 2000.

\bibitem{Johns-2007-CEK}
C.~J. Johns, J.~Mandel, A two-stage ensemble {K}alman filter for smooth data
  assimilation, {E}nvironmental and Ecological Statistics, in print, published
  online, DOI:10.1007/s10651-007-0033-0 (2007).

\bibitem{Kalman-1960-NAL}
R.~E. Kalman, A new approach to linear filtering and prediction problems,
  Transactions of the ASME -- Journal of Basic Engineering, Series D 82 (1960)
  35--45.

\bibitem{Kalnay-AMD-2003}
E.~Kalnay, Atmospheric Modeling, Data Assimilation and Predictability,
  Cambridge University Press, 2003.

\bibitem{Kolmogorov-1937-SDE}
A.~Kolmogorov, I.~Petrovskii, N.~Piscounov, A study of the diffusion equation
  with increase in the amount of substance, and its application to a biological
  problem, in: V.~M. Tikhomirov (ed.), Selected Works of {A}. {N}. {K}olmogorov
  I, Kluwer, 1991, pp. 248--270, translated by V. M. Volosov from Bull. Moscow
  Univ., Math. Mech. 1, 1--25, 1937.

\bibitem{Kremens-2003-MTT}
R.~Kremens, J.~Faulring, C.~C. Hardy, Measurement of the time-temperature and
  emissivity history of the burn scar for remote sensing applications, Paper
  J1G.5, Proceedings of the 2nd Fire Ecology Congress, Orlando FL, American
  Meteorological Society (2003).

\bibitem{Larini-1998-MFF}
M.~Larini, F.~Giroud, B.~Porterie, J.~C. Loraud, A multiphase formulation for
  fire propagation in heterogeneous combustible media, International Journal of
  Heat and Mass Transfer 41 (1998) 881--897.

\bibitem{Law-1993-CAE}
C.~K. Law, B.~H. Chao, A.~Umemura, On closure in activation energy asymptotics
  of premixed flames, Combust. Sci. Technol. 88 (1993) 59--–88.

\bibitem{Liao-2006-FCA}
W.~Liao, J.~Zhu, A.~Q.~M. Khaliq, A fourth-order compact algorithm for
  nonlinear reaction-diffusion equations with {N}eumann boundary conditions,
  Numer. Methods Partial Differential Equations 22~(3) (2006) 600--616.

\bibitem{Linn-2002-SWB}
R.~Linn, J.~Reisner, J.~J. Colman, J.~Winterkamp, Studying wildfire behavior
  using {FIRETEC}, Int. J. of Wildland Fire 11 (2002) 233--246.

\bibitem{Linn-1996-FTM}
R.~R. Linn, Transport model for prediction of wildfire behavior, ph.D. Thesis,
  Department of Mechanical Engineering, New Mexico State University (1997).

\bibitem{Mandel-2006-EIE}
J.~Mandel, Efficient implementation of the ensemble {K}alman filter, CCM Report
  231, University of Colorado Denver (2006).
\newline\urlprefix\url{http://www.math.cudenver.edu/ccm/reports/rep231.pdf}

\bibitem{Mandel-2006-PEF}
J.~Mandel, J.~D. Beezley, Predictor-corrector ensemble filters for the
  assimilation of sparse data into high dimensional nonlinear systems, CCM
  Report 232, University of Colorado Denver (2006).
\newline\urlprefix\url{http://www.math.cudenver.edu/ccm/reports/rep232.pdf}

\bibitem{Mandel-2007-DDD}
J.~Mandel, J.~D. Beezley, L.~S. Bennethum, J.~L.~C. Soham~Chakraborty, C.~C.
  Douglas, J.~Hatcher, M.~Kim, A.~Vodacek, A dynamic data driven wildland fire
  model, in: Y.~Shi, G.~D. van Albada, P.~M.~A. Sloot, J.~J. Dongarra (eds.),
  Computational Science-ICCS 2007: 7th International Conference, vol. 4487 of
  Lecture Notes in Computer Science, Springer, 2007, pp. 1042--1049.

\bibitem{Mandel-2005-TDD}
J.~Mandel, L.~S. Bennethum, M.~Chen, J.~L. Coen, C.~C. Douglas, L.~P. Franca,
  C.~J. Johns, M.~Kim, A.~V. Knyazev, R.~Kremens, V.~Kulkarni, G.~Qin,
  A.~Vodacek, J.~Wu, W.~Zhao, A.~Zornes, Towards a dynamic data driven
  application system for wildfire simulation, in: V.~S. Sunderam, G.~D. van
  Albada, P.~M.~A. Sloot, J.~J. Dongarra (eds.), Computational Science - ICCS
  2005, vol. 3515 of Lecture Notes in Computer Science, Springer, 2005, pp.
  632--639.

\bibitem{Mandel-2004-NDD}
J.~Mandel, M.~Chen, L.~P. Franca, C.~Johns, A.~Puhalskii, J.~L. Coen, C.~C.
  Douglas, R.~Kremens, A.~Vodacek, W.~Zhao, A note on dynamic data driven
  wildfire modeling, in: M.~Bubak, G.~D. van Albada, P.~M.~A. Sloot, J.~J.
  Dongarra (eds.), Computational Science - ICCS 2004, vol. 3038 of Lecture
  Notes in Computer Science, Springer, 2004, pp. 725--731.

\bibitem{Mercer-1995-CWS}
G.~N. Mercer, R.~O. Weber, Combustion wave speed, Proceedings Of the Royal
  Society Of London Series A 450 (1995) 193--198.

\bibitem{Mercer-1997-CWT}
G.~N. Mercer, R.~O. Weber, Combustion waves in two dimensions and their
  one-dimensional approximation, Combust. Theory Modelling 1 (1997) 157--165.

\bibitem{Mercer-1996-CPW}
G.~N. Mercer, R.~O. Weber, B.~F. Gray, A.~Watt, Combustion pseudo-waves in a
  system with reactant consumption and heat loss, Mathl. Comput. Modelling
  24~(8) (1996) 29--38.

\bibitem{Mickens-2005-NFD}
R.~E. Mickens, A nonstandard finite difference scheme for a {PDE} modeling
  combustion with nonlinear advection and diffusion, Math. Comput. Simulation
  69~(5-6) (2005) 439--446.

\bibitem{Morandini-2001-CRH}
F.~Morandini, P.~A. Santoni, J.~H. Balbi, The contribution of radiant heat
  transfer to laboratory-scale fire spread under the influences of wind and
  slope, Fire Safety Journal 36 (2001) 519--543.

\bibitem{Morvan-2002-BWF}
D.~Morvan, M.~Larini, J.~L. D.~P. Fernandes, A.~I. Miranda, J.~Andre,
  O.~Sero-Guillaume, D.~Calogine, P.~Cuinas, Behaviour modelling of wildland
  fires: a state of the art, {Euro-Mediterranean Wildland Fire Laboratory}, a
  `wall-less' Laboratory for Wildland Fire Sciences and Technologies in the
  Euro-Mediterranean Region (2002).

\bibitem{Norbury-1988-TCP}
J.~Norbury, A.~M. Stuart, Travelling combustion waves in a porous medium. {Part
  I}-existence, SIAM Journal on Applied Mathematics 48~(1) (1988) 155--169.

\bibitem{Norbury-1988-TCW}
J.~Norbury, A.~M. Stuart, Travelling combustion waves in a porous medium. {Part
  II}-stability, SIAM Journal on Applied Mathematics 48~(2) (1988) 374--392.

\bibitem{Ononye-2005-IFT}
A.~Ononye, A.~Vodacek, R.~Kremens, Improved fire temperature estimation using
  constrained spectral unmixing, Remote Sensing for Field Users, Proc. 10th
  Biennial USDA Forest Service Remote Sensing Applications Conference. Salt
  Lake City, UT., Am. Soc. Photogram. Remote Sens., CD-ROM (2005).

\bibitem{Ononye-2007-AEF}
A.~E. Ononye, A.~Vodacek, E.~Saber, Automated extraction of fire line
  parameters from multispectral infrared images, Remote Sensing of Environment
  108 (2007) 179--188.

\bibitem{Pastor-2003-MMC}
E.~Pastor, L.~Zarate, E.~Planas, J.~Arnaldos, Mathematical models and
  calculations systems for the study of wildland fire behavior, Prog. Energy.
  Combust. Sci. 29 (2003) 139--153.

\bibitem{Quintiere-1998-PFB}
J.~G. Quintiere, Principles of Fire Behavior, Delmar Publishers, Albany, NY,
  1998.

\bibitem{Rastigejev-2006-NSF}
Y.~Rastigejev, M.~Matalon, Numerical simulation of flames as gas-dynamic
  discontinuities, Combustion Theory and Modelling 10 (2006) 459--481.

\bibitem{Richards-1995-GMF}
G.~D. Richards, A general mathematical framework for modelling two-dimensional
  wildland fire spread, Int. J. Wildland Fire 5 (1995) 63--72.

\bibitem{Richards-1999-MMC}
G.~D. Richards, The mathematical modelling and computer simulation of wildland
  fire perimeter growth over a 3-dimensional surface, International J. of
  Wildland Fire 9~(3) (1999) 213--221.

\bibitem{Robinson-2001-IDP}
J.~C. Robinson, Infinite-dimensional dynamical systems, Cambridge Texts in
  Applied Mathematics, Cambridge University Press, Cambridge, 2001.

\bibitem{Roessler-1997-NSD}
J.~Roessler, H.~H{\"u}ssner, Numerical solution of the {$(1+2)$}-dimensional
  {F}isher's equation by finite elements and the {G}alerkin method, Math.
  Comput. Modelling 25~(3) (1997) 57--67.

\bibitem{Rothermel-1972-MMP}
R.~C. Rothermel, A mathematical model for predicting fire spread in wildland
  fires, {USDA Forest Service Research Paper INT-115} (1972).

\bibitem{Ruan-1998-EMR}
F.~Ruan, D.~McLaughlin, An efficient multivariate random field generator using
  the fast {F}ourier transform, Advances in Water Resources 21 (1998) 385--399.

\bibitem{Sembera-2001-NGM}
J.~{\v{S}}embera, M.~Bene{\v{s}}, Nonlinear {G}alerkin method for
  reaction-diffusion systems admitting invariant regions, J. Comput. Appl.
  Math. 136~(1-2) (2001) 163--176.

\bibitem{Seron-2005-EWF}
F.~J. Ser{\'o}n, D.~Guti{\'e}rrez, J.~Magall{\'o}n, L.~Ferragut, M.~I. Asensio,
  The evolution of a wildland forest fire front, Visual Computer 21 (2005)
  152--169.

\bibitem{Sethian-1999-LSM}
J.~A. Sethian, Level set methods and fast marching methods, vol.~3 of Cambridge
  Monographs on Applied and Computational Mathematics, 2nd ed., Cambridge
  University Press, Cambridge, 1999.

\bibitem{Sherratt-1998-TID}
J.~A. Sherratt, On the transition from initial data to travelling waves in the
  {F}isher-{KPP} equation, Dynam. Stability Systems 13~(2) (1998) 167--174.

\bibitem{Simeoni-2001-WIF}
A.~Simeoni, P.~A. Santoni, M.~Larini, J.~H. Balbi, On the wind advection
  influence on the fire spread across a fuel bed: modelling by a semi-physical
  approach and testing with experiments, Fire Safety Journal 36 (2001)
  491--513.

\bibitem{Sussman-1994-LSA}
M.~Sussman, P.~Smereka, S.~Osher, A level set approach for computing solutions
  to incompressible two-phase flow, J. Comput. Phys. 114 (1994) 146--159.

\bibitem{Theodoropoulos-2000-CSB}
C.~Theodoropoulos, Y.~Qian, I.~Kevrekidis, Coarse stability and bifurcation
  analysis using time-steppers: {A} reaction-diffusion example, Proc. Natl.
  Acad. Sci. USA 97 (2000) 9840–--9843.

\bibitem{Tippett-2003-ESR}
M.~K. Tippett, J.~L. Anderson, C.~H. Bishop, T.~M. Hamill, J.~S. Whitaker,
  Ensemble square root filters, Monthly Weather Review 131 (2003) 1485--1490.

\bibitem{vanLeeuwen-2003-VMF}
P.~van Leeuwen, A variance-minimizing filter for large-scale applications,
  Monthly Weather Review 131~(9) (2003) 2071--2084.

\bibitem{Viegas-2005-MMF}
D.~X. Viegas, A mathematical model for forest fires blow-up, Combustion Science
  and Technology 177 (2005) 1--25.

\bibitem{Weber-1991-MFS}
R.~O. Weber, Modelling fire spread through fuel beds, Prog. Energy Combust. 17
  (1991) 67--82.

\bibitem{Weber-1991-TCW}
R.~O. Weber, Toward a comprehensive wildfire spread model, Int. J. Wildland
  Fire 1~(4) (1991) 245--248.

\bibitem{Weber-1997-CWG}
R.~O. Weber, G.~N. Mercer, H.~S. Sidhu, B.~F. Gray, Combustion waves for gases
  ({$Le=1$}) and solids ({$Le\to\infty$}), Proceedings of the Royal Society of
  London Series A 453~(1960) (1997) 1105--1118.

\bibitem{Wotton-1999-EFF}
B.~M. Wotton, R.~S. McAlpine, M.~W. Hobbs, The effect of fire front width on
  surface fire behaviour, International Journal of Wildland Fire 9 (1999)
  247--253.

\bibitem{Zeldovich-1985-MTC}
Y.~B. Zeldovich, G.~I. Barrenblatt, V.~B. Librovich, G.~M. Makhviladze, The
  Mathematical Theory of Combustion and Explosions, Consultants Bureau, New
  York, 1985.

\bibitem{Zhao-CDS-2003}
S.~Zhao, G.~W. Wei, Comparison of the discrete singular convolution and three
  other numerical schemes for solving {F}isher's equation, SIAM J. Sci. Comput.
  25~(1) (2003) 127--147.

\bibitem{Zhou-2001-ERM}
X.~Zhou, S.~Mahalingam, Evaluation of a reduced mechanism for modeling
  combustion of pyrolysis gas in wildland fire, Combustion Science and
  Technology 171 (2001) 39--70.

\end{thebibliography}

\end{document}